# CHANGE-POINT ESTIMATION UNDER ADAPTIVE SAMPLING

By Yan Lan, Moulinath Banerjee[1] and George Michailidis[2]

*Abbott Laboratories, University of Michigan and University of Michigan*

We consider the problem of locating a jump discontinuity (change-point) in a smooth parametric regression model with a bounded covariate. It is assumed that one can sample the covariate at different values and measure the corresponding responses. Budget constraints dictate that a total of $n$ such measurements can be obtained. A multistage adaptive procedure is proposed, where at each stage an estimate of the change point is obtained and new points are sampled from its *appropriately* chosen neighborhood. It is shown that such procedures accelerate the rate of convergence of the least squares estimate of the change-point. Further, the asymptotic distribution of the estimate is derived using empirical processes techniques. The latter result provides guidelines on how to choose the tuning parameters of the multistage procedure in practice. The improved efficiency of the procedure is demonstrated using real and synthetic data. This problem is primarily motivated by applications in engineering systems.

**1. Introduction.** The problem of estimating the location of a jump discontinuity (change-point) in an otherwise smooth curve has been extensively studied in the nonparametric regression and survival analysis literature [see, e.g., Dempfle and Stute (2002), Gijbels, Hall and Kneip (1999), Hall and Molchanov (2003), Kosorok and Song (2007), Koul and Qian (2002), Loader (1996), Müller (1992), Müller and Song (1997), Pons (2003), Ritov (1990) and references therein]. In the classical setting, measurements on all $n$ covariate-response pairs are available *in advance*, and the main issue is to estimate as accurately as possible the location of the change-point. However, there are applications where it is possible to sample the response at any covariate

Received August 2007; revised December 2007.
[1]Supported in part by NSF Grant DMS-07-05288.
[2]Supported in part by NIH Grant P41-18627.
*AMS 2000 subject classifications.* 62F12, 62K99.
*Key words and phrases.* Adaptive sampling, change point estimation, multistage procedure, Skorokhod topology, two-stage procedure, zoom-in.

value of the experimenter's choice. The only hard constraint is that the total budget of measurements to be obtained is fixed a priori.

For example, consider the following example from system engineering. There is a stochastic flow of jobs/customers of various types arriving to the system with random service requests. Jobs waiting to be served are placed in queues of infinite capacity. The system's resources are allocated to the various job classes (queues) according to some service policy. This system serves as a canonical queueing model for many applications, including network switches, flexible manufacturing systems, wireless communications, etc. [Hung and Michailidis (2008)]. A quantity of great interest to the system's operator is the average delay of the customers, which is a key performance metric of the quality of service offered by the system.

The average delay of the customers in a two-class system as a function of its loading, for a resource allocation policy introduced and discussed in Hung and Michailidis (2008), is shown in Figure 1. Specifically, the system was simulated under 134 loading settings and fed by input/service request processes obtained from real network traces and the average delay of 500,000 customers recorded. It can be seen that for loading around 0.8 there is a marked discontinuity in the response, which indicates that under the specified resource-allocation policy the service provided to the customers deteriorates. It is of interest to locate the "threshold" where such a change in the quality of service occurs. It should be pointed out that this threshold would occur at different system loadings for different allocation policies.

A few comments on the setting implied by this example are in order. First, the experimenter can *select* covariate values (in this case, the system's loading) and subsequently obtain their corresponding sample responses. Second, the sampled responses are *expensive* to obtain; for example, the average delay is obtained by running a fairly large-scale discrete event simulation of the system under consideration, involving half a million customers. For systems, comprising of a large number of customer classes, more computationally intensive simulations that can last days must be undertaken. Third, in many situations there is an a priori *fixed budget* of resources; for this example it may correspond to CPU time, in other engineering applications to emulation time, while in other scientific contexts to real money.

Given the potentially limited budget of points that can be sampled and lack of a priori knowledge about the location of the change-point, the following strategy looks promising. A certain portion of the budget is used to obtain an initial estimate of the change-point based on a least squares criterion. Subsequently, a neighborhood around this initial estimate is specified and the remaining portion of the available points are sampled from it, together with their responses, that yield a new estimate of the change-point. Intuition suggests that if the first-stage estimate is fairly accurate, the more intensive sampling in its neighborhood ought to produce a more accurate estimate than the one that would have been obtained by laying out the entire budget of points in a uniform fashion. Obviously, the procedure with its "zoom-in" characteristics can be extended beyond two stages.

The goal of this paper is to formally introduce such multistage adaptive procedures for change-point estimation and examine their properties. In particular, the following important issues are studied and resolved: (i) the selection of the size of the neighborhoods, (ii) the rate of convergence of the multistage least squares estimate, together with its asymptotic distribution and (iii) allocation of the available budget at each stage.

The proposed procedure should be contrasted with the well-studied sequential techniques for change-point detection, since the underlying setting exhibits marked differences. In its simplest form, the sequential change-point detection problem can be formulated as follows: there is a process that generates a sequence of independent observations $X_1, X_2, \ldots$ from some distribution $F_0$. At some unknown point in time $\tau$, the distribution changes and hence observations $X_\tau, X_{\tau+1}, \ldots$ are generated from $F_1$. The objective is to raise an alarm as soon as the data-generating mechanism switches to a new distribution. This problem originally arose in statistical quality control and over the years has found important applications in other fields.

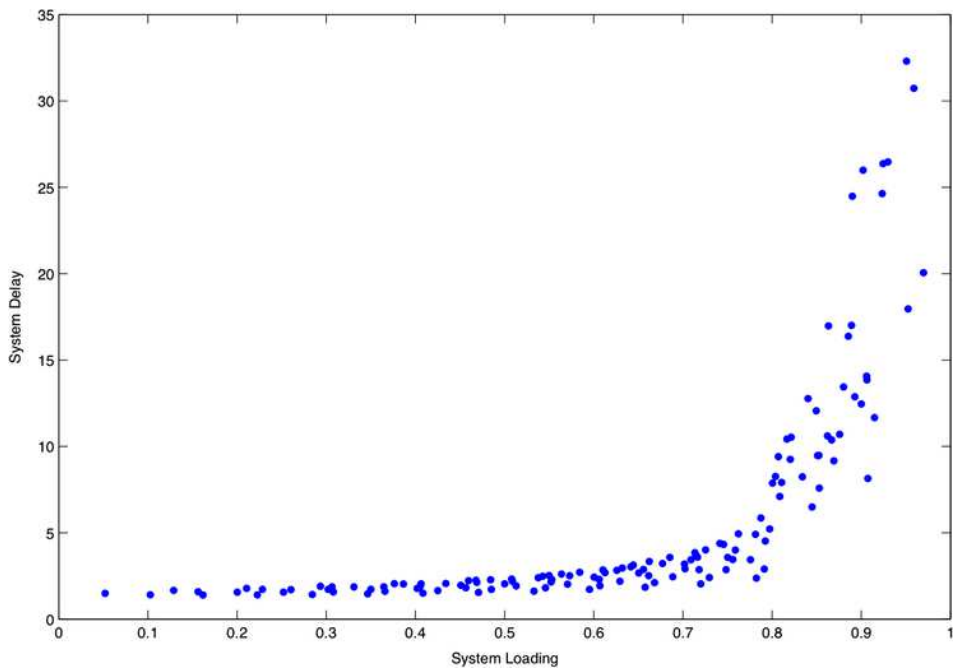

FIG. 1. *Average delay as a function of system loading for a two-class parallel processing system.*

Being a canonical problem in sequential analysis, many detection procedures have been proposed in the literature over the years in discrete and continuous time, under various assumptions on the distribution of $\tau$ and the data-generating mechanism. The literature on this subject is truly enormous; a comprehensive treatment of the problem can be found in the book by Basseville and Nikiforov (1993), while some recent developments and new challenges are discussed in the review paper by Lai (2001). An important difference in our setting is the control that the experimenter exercises over the data generation process and also the absence of physical time, a crucial element in the sequential change-point problem.

The remainder of the paper is organized as follows. In Section 2, the one-stage procedure is briefly reviewed. In Section 3, the proposed procedure based on adaptive sampling is introduced and the main results regarding the rate of convergence of the new estimator of the change point and its asymptotic distribution are established; in addition, various generalizations are discussed. In Section 4, the performance of the proposed estimator in finite samples is studied through an extensive simulation and practical guidelines are discussed for its various tuning parameters. In Section 5, various techniques for constructing confidence intervals for the change-point, based on the two-stage procedure, are discussed. Finally, the majority of the technical details are provided in the Appendix.

**2. The classical problem.** In this study, we focus on parametric models for the regression function of the type

$$Y_i = \mu(X_i) + \epsilon_i, \qquad i = 1, 2, \ldots, n,$$

where

$$(1) \qquad \mu(x) = \psi_l(\beta_l, x) 1(x \leq d^0) + \psi_u(\beta_u, x) 1(x > d^0)$$

with $\psi_l(\beta_l, x)$ and $\psi_u(\beta_u, x)$ are both (at least) twice continuously differentiable in $\beta$ and infinitely differentiable in $x$ and $\psi_l(\beta_l, d^0) \neq \psi_u(\beta_u, d^0)$, so that $d^0$ is the unique point of discontinuity—a change point—of the regression function.

The $\epsilon_i$'s are assumed to be i.i.d. symmetric mean 0 errors with common (unknown) error variance $\sigma^2$ and are independent of the $X_i$'s, which are i.i.d. and are distributed on $[0, 1]$ according to some common density $f_X(\cdot)$. The simplest possible parametric candidate for $\mu(x)$, which we will focus on largely to illustrate the key ideas in the paper, is the simple step function $\mu(x) = \alpha_0 1(x \leq d^0) + \beta_0 1(x > d^0)$.

Estimating $d^0$ based on the above data is coined as the "classical problem." A standard way to estimate the parameters $(\beta_l, \beta_u, d^0)$ is to solve a least-squares problem. We start by introducing some necessary notation.

Let $\mathbb{P}_n$ denote the empirical distribution of the data vector $\{X_i, Y_i\}_{i=1}^n$ and $P$ the true distribution of $(X_1, Y_1)$. For a function $f$ defined on the space $[0,1] \times \mathbb{R}$ (in which the vector $(X_1, Y_1)$ assumes values) and a measure $Q$ defined on the Borel $\sigma$-field on $[0,1] \times \mathbb{R}$, we denote $\int f \, dQ$ as $Qf$. We now turn our attention to the least-squares problem.

The objective is to minimize $\mathbb{P}_n[(y - \psi_l(\alpha, x))^2 1(x \leq d) + (y - \psi_u(\beta_u, x))^2 1 \times (x > d)]$ over all $(\alpha, \beta, d)$, with $0 \leq d \leq 1$. Let $(\tilde{\beta}_{l,n}, \tilde{\beta}_{u,n}, \tilde{d}_n)$ denote a vector of minimizers. Note that we refer to "a vector" of minimizers, since there will be, in general, multiple tri-vectors that minimize the criterion function. The asymptotic properties of such a vector can be studied by invoking either the methods of Pons (2003) or those (in Chapter 14) of Kosorok (2008) and Kosorok and Song (2007). We do not provide the details, but state the results that are essential to the *multistage learning procedures* that we formulate in the next section. We clarify next the meaning of a minimizer of a right-continuous, real-valued function with left limits (say $f$) defined on an interval $I$. Specifically, any point $z \in I$ that satisfies $f(z) \wedge f(z-) = \min_{w \in I} f(w)$ is defined to be a minimizer of $f$. Also, in order to discuss the results for the classical procedure and those for the proposed multistage procedures, we need to define a family of compound Poisson processes that arise in the description of the asymptotic properties of the estimators of the change point.

*A family of compound Poisson processes.* For a positive constant $\Lambda$, let $\nu^+(\cdot)$ be a Poisson process on $[0, \infty)$ with right-continuous sample paths, with $\nu^+(s) \sim \text{Poi}(\Lambda s)$ for $s > 0$. Let $\tilde{\nu}^+(\cdot)$ be another independent Poisson process on $[0, \infty)$ with left-continuous sample paths, with $\tilde{\nu}^+(s) \sim \text{Poi}(\Lambda s)$ and define a (right-continuous) Poisson process on $(-\infty, 0]$ by $\{\nu^-(s) = \tilde{\nu}^+(-s) : s \in (-\infty, 0]\}$. Let $\{\eta_i\}_{i=1}^\infty$ and $\{\eta_{-i}\}_{i=1}^\infty$ be two independent sequences of i.i.d. random variables where each $\eta_j$ ($j$ assumes both positive and negative values) is distributed like $\eta$, $\eta$ being a mean 0 random variable with finite variance $\rho^2$. Given a positive constant $A$, define families of random variables $\{V_i^+\}_{i=1}^\infty$ and $\{V_i^-\}_{i=1}^\infty$ where, for each $i \geq 1$, $V_i^+ = A/2 + \eta_i$ and $V_i^- = -A/2 + \eta_{-i}$. Set $V_0^+ = V_0^- \equiv 0$. Next, define compound Poisson processes $\mathbb{M}_1$ and $\mathbb{M}_2$ on $(-\infty, \infty)$ as follows: $\mathbb{M}_1(s) = (\sum_{0 \leq i \leq \nu^+(s)} V_i^+) 1(s \geq 0)$ and $\mathbb{M}_2(s) = (\sum_{0 \leq i \leq \nu^-(s)} V_i^-) 1(s \leq 0)$. Finally, define the two-sided compound Poisson process $\mathbb{M}_{A,\eta,\Lambda}(s) = \mathbb{M}_1(s) - \mathbb{M}_2(s)$. It is not difficult to see that $\mathbb{M}$, almost surely, has a minimizer (in which case it has multiple minimizers, since the sample paths are piecewise constant). Let $d_l(A, \eta, \Lambda)$ denote the smallest minimizer of $\mathbb{M}_{A,\eta,\Lambda}$ and $d_u(A, \eta, \Lambda)$ its largest one, which are almost surely well defined. Then, the following relation holds:

$$(2) \qquad (d_l(A, \eta, \Lambda), d_u(A, \eta, \Lambda)) \stackrel{d}{=} \left( d_l\left(\frac{A}{\rho}, \frac{\eta}{\rho}, \Lambda\right), d_u\left(\frac{A}{\rho}, \frac{\eta}{\rho}, \Lambda\right) \right)$$

$$(3) \qquad \stackrel{d}{=} \frac{1}{\Lambda} \left( d_l\left(\frac{A}{\rho}, \frac{\eta}{\rho}, 1\right), d_u\left(\frac{A}{\rho}, \frac{\eta}{\rho}, 1\right) \right).$$

For the "classical problem," the following proposition holds.

PROPOSITION 1. *Consider the model described at the beginning of Section 2. Suppose that $X$ has a positive bounded density on $[0,1]$ and that $d^0$ is known to lie in the interval $[\epsilon_0, 1-\epsilon_0]$ for some small $\epsilon_0 > 0$. Let $(\hat{\beta}_{l,n}, \hat{\beta}_{u,n}, \hat{d}_n)$ denote that minimizing tri-vector $(\tilde{\beta}_{l,n}, \tilde{\beta}_{u,n}, \tilde{d}_n)$, for which the third component is minimal. Then, $(\sqrt{n}(\hat{\beta}_{l,n} - \beta_l), \sqrt{n}(\hat{\beta}_{u,n} - \beta_u), n(\hat{d}_n - d^0))$ is $O_p(1)$. Furthermore, the first two components of this vector are asymptotically independent of the third and*

$$n(\hat{d}_n - d^0) \xrightarrow{d} d_l(|\mu(d^0+) - \mu(d^0)|, \epsilon_1, f_X(d^0))$$

$$\stackrel{d}{\equiv} \frac{1}{f_X(d^0)} d_l\left(\frac{|\mu(d^0+) - \mu(d^0)|}{\sigma}, \frac{\epsilon_1}{\sigma}, 1\right).$$

*Heteroscedastic errors.* The proposition can be generalized readily to cover the case of heteroscedastic errors. A generalization of the classical model to the heteroscedastic case is as follows: we observe $n$ i.i.d. observations from the model $Y = \mu(X) + \sigma(X)\tilde{\epsilon}$ where $\mu(x)$ is as defined in (1), $\tilde{\epsilon}$ and $X$ are independent, $\tilde{\epsilon}$ is a mean 0 error with unit variance and $\sigma^2(x)$ is a variance function (assumed continuous). As in the homoscedastic case, an unweighted least-squares procedure is used to estimate the parameters $(\beta_l, \beta_u, d^0)$. As before, letting $(\hat{\beta}_{l,n}, \hat{\beta}_{u,n}, \hat{d}_n)$ denote that minimizing tri-vector $(\tilde{\beta}_{l,n}, \tilde{\beta}_{u,n}, \tilde{d}_n)$, for which the third component is minimal, we have:

$$n(\hat{d}_n - d^0) \xrightarrow{d} d_l(|\mu(d^0+) - \mu(d^0)|, \sigma(d^0)\tilde{\epsilon}, f_X(d^0))$$

$$\stackrel{d}{\equiv} \frac{1}{f_X(d^0)} d_l\left(\frac{|\mu(d^0+) - \mu(d^0)|}{\sigma(d^0)}, \tilde{\epsilon}, 1\right).$$

**3. The two-stage procedure.** We first describe a two-stage procedure for estimating the (unique) change-point. In what follows, we consider a regression scenario where the response $Y_x$ generated at covariate level x can be written as $Y_x = \mu(x) + \epsilon$, where $\epsilon$ is a error variable with finite variance and $\mu$ is the regression function. The errors corresponding to different covariate levels are i.i.d. We first focus on the simple regression function $\mu(x) = \alpha_0 1(x \leq d^0) + \beta_0 1(x > d^0)$ and discuss generalizations to more complex parametric models later on. We are allowed to sample $n$ covariate-response pairs at most and are free to sample a response from any covariate level that we like.

- Step 1. At stage one, $\lambda n$ covariate values are sampled uniformly from $[0, 1]$ and responses are obtained. Denote the observed data by $\{X_i, Y_i\}_{i=1}^{n_1}$, $n_1 = \lambda n$ and the corresponding estimated location of the change-point by $\hat{d}_{n_1}$.

- Step 2. Sample the remaining $n_2 = (1-\lambda)n$ covariate-response pairs $\{U_i, W_i\}_{i=1}^{n_2}$, where:

$$W_i = \mu(U_i) + \epsilon_i, \qquad U_i \sim \text{Unif}[\hat{a}_{n_1}, \hat{b}_{n_1}]$$

and $[\hat{a}_{n_1}, \hat{b}_{n_1}] = [\hat{d}_{n_1} - K n_1^{-\gamma}, \hat{d}_{n_1} + K n_1^{-\gamma}]$, $0 < \gamma < 1$ and $K$ is some constant.

Obtain an updated estimate of the change-point based on the $n_2$ covariate-response pairs from stage two, which is denoted by $\hat{d}_{n_2}$.

We discuss the basic procedure in some more detail. Let $(\hat{\alpha}_{n_1}, \hat{\beta}_{n_1}, \hat{d}_{n_1})$ denote the parameter estimates obtained from stage one. Let $\mathbb{P}_{n_2}$ denote the empirical measure of the data points $\{U_i, W_i\}_{i=1}^{n_2}$. The updated estimates are computed by minimizing

$$\mathbb{P}_{n_2}[\{(w - \hat{\alpha}_{n_1})^2 I(u \leq d) + (w - \hat{\beta}_{n_1})^2 I(u > d)],$$

which, as is readily seen, is equivalent to minimizing the process

$$\tilde{\mathbb{M}}_{n_2}(d) \equiv \mathbb{P}_{n_2}[\{(w - \hat{\alpha}_{n_1})^2 - (w - \hat{\beta}_{n_1})^2\}(I(u \leq d) - I(u \leq d^0))].$$

The process $\tilde{\mathbb{M}}_{n_2}$ is a piecewise-constant, right-continuous function with left limits. We let $\hat{d}_{n_2,l}$ and $\hat{d}_{n_2,u}$ denote its minimal and maximal minimizers, respectively. Our goal is to determine the joint limit distribution of normalized versions of $(\hat{d}_{n_2,l}, \hat{d}_{n_2,u})$. This is described in the theorems that follow.

THEOREM 1. *Assume that the error variable $\epsilon$ in the regression model has a finite moment generating function in a neighborhood of 0. Then, the random vector $n^{1+\gamma}(\hat{d}_{n_2,l} - d^0, \hat{d}_{n_2,u} - d^0)$ is $O_p(1)$.*

REMARK. The proof of this theorem is fairly technical and particularly long and thus deferred to the Appendix. However, a few words regarding the intuition behind the *accelerated rate* of convergence are in order. For simplicity, consider sampling procedures where instead of sampling from a uniform distribution on the interval of interest, sampling takes place on a uniform grid on the interval. The interval from which sampling takes place at the second stage has length $2K n_1^{-\gamma}$. Since the $n_2$ covariate values are equispaced over this interval, the resolution of the resulting grid at which responses are measured is $O(n_1^{-\gamma}/n_2) = O(n^{-(1+\gamma)})$ and this determines the rate of convergence of the two-stage estimator (just as the rate of convergence in the classical procedure where $n$ covariates are equispaced over $[0,1]$ is given by the resolution of the resulting grid in that situation, which is simply $(n^{-1})$).

We next describe the limit distributions of the normalized estimates considered in Theorem 1.

THEOREM 2. *Set $C(K, \lambda, \gamma) = (2K)^{-1}(\lambda/(1-\lambda))^\gamma$. The random vector $n_2^{1+\gamma}(\hat{d}_{n_2,l} - d^0, \hat{d}_{n_2,u} - d^0)$ converges in distribution to*

$$(d_l(|\alpha_0 - \beta_0|, \epsilon, C(K, \lambda, \gamma)), d_u(|\alpha_0 - \beta_0|, \epsilon, C(K, \lambda, \gamma))).$$

REMARK. The asymptotic distributions of the "zoom-in" estimators are given by the minimizers of a compound Poisson process. The underlying Poisson process is basically the limiting version of the count process $\{\mathcal{P}_n(s) : s \in \mathbb{R}\}$, where $\mathcal{P}_n(s)$ counts the number of $U_i$'s in the interval $(d^0, d^0 + s/n_2^{1+\gamma}] \cup (d^0 + s/n_2^{1+\gamma}, d^0]$. It can be readily checked that marginally, $\mathcal{P}_n(s)$, converges in distribution to a Poisson random variable with mean $C(K, \lambda, \gamma)s$, using the Poisson approximation to the Binomial distribution. On the other hand, the size of the jumps of the compound Poisson process is basically determined by $|\alpha_0 - \beta_0|/\sigma$, the signal-to-noise ratio in the model.

*General parametric models.* These results admit ready extensions to the case where the function $\mu(x)$ is as defined in (1). As in the case of a piecewise constant $\mu$, $n_1 \equiv \lambda n$ points are initially used to obtain least-squares estimates of $(\beta_l, \beta_u, d^0)$, which we denote by $(\hat{\beta}_{l,n_1}, \hat{\beta}_{u,n_1}, \hat{d}_{n_1})$. Step 2 of the two-stage procedure is identical and the updated estimate $\hat{d}_{n_2}$ is computed by minimizing the criterion function

$$\mathbb{P}_{n_2}[(w - \psi_l(\hat{\beta}_{l,n_1}, u))^2 I(u \leq d) + (w - \psi_u(\hat{\beta}_{u,n_1}, u))^2 I(u > d)],$$

which is equivalent to minimizing

$$\tilde{\mathbb{M}}_{n_2}(d) = \mathbb{P}_{n_2}[\{(w - \psi_l(\hat{\beta}_{l,n_1}, u))^2 - (w - \psi_u(\hat{\beta}_{u,n_1}, u))^2\} \times (I(u \leq d) - I(u \leq d^0))].$$

Letting $\hat{d}_{n_2,l}$ and $\hat{d}_{n_2,u}$ denote the smallest and largest argmins of $\tilde{\mathbb{M}}_{n_2}$, respectively (as in the piecewise-constant function case), we have the following proposition.

PROPOSITION 2. *The random vector $n_2^{1+\gamma}(\hat{d}_{n_2,l} - d^0, \hat{d}_{n_2,u} - d^0)$ converges in distribution to*

$$(d_l(|\psi_l(\beta_l, d^0) - \psi_u(\beta_u, d^0)|, \epsilon, C(K, \lambda, \gamma)),$$
$$d_u(|\psi_l(\beta_l, d^0) - \psi_u(\beta_u, d^0)|, \epsilon, C(K, \lambda, \gamma))).$$

*The heteroscedastic case.* Similar results continue to hold for a heteroscedastic regression model. We formulate the heteroscedastic setting as follows. At any given covariate level $x$, the observed response $Y_x = \mu(x) + \sigma(x)\tilde{\epsilon}$ with $\mu(x)$ as defined in (1), $\sigma^2(x)$ is a (continuous) variance function and $\tilde{\epsilon}$ is a symmetric error variable with unit variance. The errors corresponding to different covariate values are independent. Using the same two-stage procedure as described above, the following proposition obtains.

PROPOSITION 3. *We have*
$$n_2^{1+\gamma}(\hat{d}_{n_2,l} - d^0, \hat{d}_{n_2,u} - d^0)$$
$$\xrightarrow{d} (d_l(|\psi_l(\beta_l, d^0) - \psi_u(\beta_u, d^0)|, \sigma(d^0)\tilde{\epsilon}, C(K, \lambda, \gamma)),$$
$$d_u(|\psi_l(\beta_l, d^0) - \psi_u(\beta_u, d^0)|, \sigma(d^0)\tilde{\epsilon}, C(K, \lambda, \gamma))).$$

REMARK. With choice of a constant variance function, $\sigma^2(x) \equiv \sigma^2$, the heteroscedastic model reduces to the homoscedastic one. We nevertheless present results for these two situations separately. We also subsequently derive our results for the homoscedastic case, the derivations extending almost trivially to the heteroscedastic case.

3.1. *Proof of Theorem 2.* For the proof of this theorem (and the proof of Lemma 3.2 in the Appendix) we denote the process $\mathbb{M}_{|\alpha_0-\beta_0|,\epsilon,C(K,\lambda,\gamma)}$ simply by $\mathbb{M}$ and its smallest and largest minimizers simply by $(d_l, d_u)$. Our proof of this theorem will rely on continuous mapping for the argmin functional. For the sake of concreteness, in what follows, we assume that $\alpha_0 < \beta_0$. Under this assumption, with probability increasing to 1 as $n$ (and consequently $n_1$) goes to infinity, $\hat{\alpha}_{n_1} < \hat{\beta}_{n_1}$ and $d^0$ belongs to the set $[\hat{d}_{n_1} - Kn_1^{-\gamma}, \hat{d}_{n_1} + Kn_1^{-\gamma}]$. On this set $(\hat{d}_{n_2,l}, \hat{d}_{n_2,u})$ can be obtained by minimizing (the equivalent) criterion function

$$\mathbb{P}_{n_2}\left[\left(w - \frac{\hat{\alpha}_{n_1} + \hat{\beta}_{n_1}}{2}\right)(I(u \leq d) - I(u \leq d^0))\right]$$

and $d^0$ is characterized as:

$$d^0 = \arg\min \mathbf{P}\left[\left(w - \frac{\hat{\alpha}_{n_1} + \hat{\beta}_{n_1}}{2}\right)(I(u \leq d) - I(u \leq d^0))\right],$$

where $\mathbf{P}$ is the distribution of $(W, U)$. Therefore, in what follows, we take:

$$\tilde{\mathbb{M}}_{n_2}(d) = \mathbb{P}_{n_2}\left[\left(w - \frac{\hat{\alpha}_{n_1} + \hat{\beta}_{n_1}}{2}\right)(I(u \leq d) - I(u \leq d^0))\right]$$

and $\hat{d}_{n_2,l}$ and $\hat{d}_{n_2,u}$ to be the smallest and largest argmins of this stochastic process. Set $(\xi_{n,l}, \xi_{n,u}) = n_2^{1+\gamma}(\hat{d}_{n_2,l} - d^0, \hat{d}_{n_2,u} - d^0)$. Then $(\xi_{n,l}, \xi_{n,u})$ is the vector of smallest and largest argmins of the stochastic process:

$$\mathbb{M}_{n_2}(s) = \sum_{i=1}^{n_2}\left[\left(W_i - \frac{\hat{\alpha}_{n_1} + \hat{\beta}_{n_1}}{2}\right)\left(I\left(U_i \leq d^0 + \frac{s}{n_2^{1+\gamma}}\right) - I(U_i \leq d^0)\right)\right]$$
$$= \mathbb{M}_{n_2}^+(s) - \mathbb{M}_{n_2}^-(s),$$

where $\mathbb{M}_{n_2}^+(s) = \mathbb{M}_{n_2}(s)1(s \geq 0)$ and $\mathbb{M}_{n_2}^-(s) = -\mathbb{M}_{n_2}(s)1(s \leq 0)$.

We now introduce some notation that is crucial to the subsequent development. Let $\mathcal{S}$ denote the class of piecewise-constant, right-continuous functions with left limits (from $\mathbb{R}$ to $\mathbb{R}$) that are continuous at every integer point, assume the value 0 at 0 and possess finitely many jumps in every compact interval $[-C, C]$ where $C > 0$ is an integer. Let $\tilde{f}$ denote the pure jump process (of jump size 1) corresponding to the function $f$; that is, $\tilde{f}$ is the piecewise-constant, right-continuous function with left limits, such that for any $s > 0$, $\tilde{f}(s)$ counts the number of jumps of the function $f$ in the interval $[0, s]$, while for $s < 0$, $\tilde{f}(s)$ counts the number of jumps in the set $(s, 0)$.

For any positive integer $C > 0$, let $\mathcal{D}[-C, C]$ denote the class of all right-continuous functions with left limits with domain $[-C, C]$ equipped with the Skorokhod topology and let $\mathcal{D}[-C, C] \times [-C, C]$ denote the corresponding product space. Finally, let $\mathcal{D}_C^0$ denote the (metric) subspace of $([-C, C]) \times ([-C, C])$ that comprises all function pairs of the form $(f|_{[-C,C]}, \tilde{f}|_{[-C,C]})$ for $f \in \mathcal{S}$. We have the following lemma that is proved in the Appendix section of Lan, Banerjee and Michailidis (2007).

LEMMA 3.1. *Let $\{f_n\}$ and $f_0$ be functions in $\mathcal{S}$, such that for every positive integer $C$, $(f_n|_{[-C,C]}, \tilde{f}_n|_{[-C,C]})$ converges to $(f_0|_{[-C,C]}, \tilde{f}_0|_{[-C,C]})$ in $\mathcal{D}_0^C$ where $f_0$ satisfies the property that no two flat stretches of $f_0$ have the same height. Let $l_{n,C}$ and $u_{n,C}$ denote the smallest and the largest minimizers of $f_n$ on $[-C, C]$, and $l_{0,C}$ and $u_{0,C}$ denote the corresponding functionals for $f_0$. Then $(l_{n,C}, u_{n,C}) \to (l_{0,C}, u_{0,C})$.*

Consider the sequence of stochastic processes $\mathbb{M}_{n_2}(s)$ and let $\mathbb{J}_{n_2}(s)$ denote the corresponding jump processes. We have:

$$\mathbb{J}_{n_2}(s) = \text{sign}(s) \sum_{i=1}^{n_2} \left[ \left( I\left(U_i \leq d^0 + \frac{s}{n_2^{1+\gamma}}\right) - I(U_i \leq d^0) \right) \right]$$
$$= \mathbb{J}_{n_2}^+(s) + \mathbb{J}_{n_2}^-(s),$$

where $\mathbb{J}_{n_2}^+(s) = \mathbb{J}_{n_2}(s) 1(s \geq 0)$ and $\mathbb{J}_{n_2}^-(s) = \mathbb{J}_{n_2}(s) 1(s \leq 0)$. The jump process corresponding to $\mathbb{M}(s)$ is denoted by $\mathbb{J}(s)$ and is given by $\nu^+(s) 1(h \geq 0) + \nu^-(s) 1(h \leq 0)$. For each $n$, $\{\mathbb{M}_{n_2}(s) : s \in \mathbb{R}\}$ lives in $\mathcal{S}$ with probability one. Also, with probability one, $\{\mathbb{M}(s) : s \in \mathbb{R}\}$ lives in $\mathcal{S}$. Also, on a set of probability one (which does not depend on $C$), for every positive integer $C$, $((\mathbb{M}_{n_2}(s), \mathbb{J}_{n_2}(s)) : s \in [-C, C])$ belongs to $\mathcal{D}_0^C$ and so does $((\mathbb{M}(s), \mathbb{J}(s)) : s \in [-C, C])$. Let $(\xi_{n,C,l}, \xi_{n,C,u})$ denote the smallest and largest arg min of $\mathbb{M}_{n_2}$ restricted to $[-C, C]$ and let $(d_{C,l}, d_{C,u})$ denote the corresponding functionals for $\mathbb{M}$ restricted to $[-C, C]$. We prove in the Appendix:

LEMMA 3.2. *For every $C > 0$, $((\mathbb{M}_{n_2}(s), \mathbb{J}_{n_2}(s)) : s \in [-C, C])$ converges in distribution to $((\mathbb{M}(s), \mathbb{J}(s)) : s \in [-C, C])$ in the space $\mathcal{D}_C^0$.*

Consider the function $h$ that maps an element (a pair of functions) of $\mathcal{D}_0^C$ to the two-dimensional vector given by the smallest arg min and the largest arg min of the first component of the element. Using the fact that almost surely no two flat stretches of $\mathbb{M}$ have the same height, it follows by Lemma 3.1 that the process $((\mathbb{M}(s), \mathbb{J}(s)) : s \in [-C, C])$ belongs, almost surely, to the continuity set of the function $h$. This, coupled with the distributional convergence established in Lemma 3.2 leads to the conclusion that

$$(4) \qquad (\xi_{n,C,l}, \xi_{n,C,u}) \xrightarrow{d} (d_{C,l}, d_{C,u}).$$

We will show that $(\xi_{n,l}, \xi_{n,u}) \to (d_l, d_u)$. To this end, we use the following lemma.

LEMMA 3.3. *Suppose that $\{W_{n\epsilon}\}, \{W_n\}$ and $\{W_\epsilon\}$ are three sets of random vectors such that:*
  (i) $\lim_{\epsilon \to 0} \limsup_{n \to \infty} P[W_{n\epsilon} \neq W_n] = 0$,
  (ii) $\lim_{\epsilon \to 0} P[W_\epsilon \neq W] = 0$ *and*
  (iii) *for every $\epsilon > 0$, $W_{n\epsilon} \xrightarrow{d} W_\epsilon$ as $n \to \infty$.*
*Then, $W_n \xrightarrow{d} W$, as $n \to \infty$.*

Before applying the lemma, we first note the following facts: (a) The sequence of (smallest and largest minimizers) $(\xi_{n,l}, \xi_{n,u})$ is $O_p(1)$, and (b) The minimizers $(d_l, d_u)$ are $O_p(1)$. Now, in the above lemma, set $\epsilon = 1/C$, $W_{n\epsilon} = (\xi_{n,C,l}, \xi_{n,C,u})$, $W_\epsilon = (d_{C,l}, d_{C,u})$, $W_n = (\xi_{n,l}, \xi_{n,u})$ and $W = (d_l, d_u)$. Condition (iii) is established in (4). From (a) and (b) it follows that conditions (i) and (ii) of the lemma are satisfied. We conclude that $(\xi_{n,l}, \xi_{n,u}) \xrightarrow{d} (d_l, d_u)$.
□

REMARK. It is instructive to compare the obtained result on the convergence of the (nonunique) arg min functional to that considered in Ferger (2004). Ferger deals with the convergence of the arg max functional under the Skorokhod topology in Theorems 2 and 3 of his paper. Since the arg max functional is *not continuous* under the Skorokhod topology, an exact result on distributional convergence cannot be achieved. Instead, asymptotic upper and lower bounds are obtained on the distribution function of the arg max in terms of the smallest maximizer and the largest maximizer of the limit process [page 88 of Ferger (2004)]. The result we obtain here is, admittedly, in a more specialized set-up than the one considered in his paper, but it is stronger since we are able to show *exact* distributional convergence of argmins. This is achieved at the cost of some extra effort: establishing the joint convergence of the original processes, whose argmins are of interest, and their jump processes, and subsequently invoking continuous mapping. Under this stronger mode of convergence, the arg min functional indeed turns out to be continuous, as Lemma 3.1 shows [the arguments employed are similar in spirit to those in Section 14.5.1 of Kosorok (2008)]. This result allows us to construct asymptotic confidence intervals that have *exact* coverage at any given level, as opposed to the conservative intervals proposed in Ferger (2004). That the exact confidence intervals buy us significant precision over the conservative ones is evident from the reported simulation results discussed in Section 5.

## 4. Multistage procedures and strategies for parameter allocation.

4.1. *Multistage procedures.* We consider a generalization of the two-stage procedure to $P$ stages in the setting of the heteroscedastic model with a general parametric regression function $\mu$. Let $\lambda_1, \lambda_2, \ldots, \lambda_P$ be the proportions of points used at each stage (where $\lambda_1 + \lambda_2 + \cdots + \lambda_P = 1$) and let $n_i = \lambda_i n$. Also, fix sequences of numbers $0 < \gamma_{P-1} < \cdots < \gamma_1 < 1$ and $K_1, K_2, \ldots, K_{P-1}$ (with $K_i > 0$). Having used $n_1$ points to construct the initial estimate $\hat{d}_{n_1}$, in the $q$th ($2 \leq q \leq P$) stage, define the sampling neighborhood as $[\hat{d}_{n_{q-1}} - K_{q-1} n_{q-1}^{-((q-2)+\gamma_{q-1})}, \hat{d}_{n_{q-1}} + K_{q-1} n_{q-1}^{-((q-2)+\gamma_{q-1})}]$, sample $n_q$ covariate-response pairs $\{W_i, U_i\}_{i=1}^{n_q}$ from this neighborhood: $W_i = \mu(U_i) + \epsilon_i$ and update the estimate of the change-point to $\hat{d}_{n_q}$. Let $(\hat{d}_{n_P,l}, \hat{d}_{n_P,u})$ denote the smallest and largest estimates at stage $P$. It can be shown using analogous arguments as those in Theorems 3.1 and 3.2 that

$$n_P^{(P-1)+\gamma_{P-1}}((\hat{d}_{n_P,l} - d^0), (\hat{d}_{n_P,u} - d^0))$$

is $O_p(1)$ and converges in distribution to $(d_l, d_u)$, where $(d_l, d_u)$ is the vector of the smallest and the largest argmins of the process

$$\mathbb{M}(|\psi_l(\beta_l, d^0) - \psi_u(\beta_u, d^0)|, \sigma(d^0)\tilde{\epsilon}, C_P)$$

with $C_P = (1/2 K_{P-1})(\lambda_{P-1}/\lambda_P)^{((P-2)+\gamma_{P-1})}$.

4.2. *Strategies for parameter allocation.* In this section, we describe strategies for selecting the tuning parameters $K, \gamma$ and $\lambda$ used in the procedure. We do this in the setting of the simple regression model $\mu(x) = \alpha_0 1(x \leq d^0) + \beta_0 1(x > d^0)$ and homoscedastic normal errors, obvious analogues holding in more general settings.

Recall that $(\hat{d}_{n_q,l}, \hat{d}_{n_q,u})$ are the minimal and maximal minimizers at step $q, 2 \leq q \leq P$. Set: $\hat{d}_{q,av} = (\hat{d}_{n_q,l} + \hat{d}_{n_q,u})/2$. In what follows we use this as our $q$th stage estimate of the change-point.

We start with the case of $P = 2$ to fix ideas and motivate the general result. Using notation from Theorem 2 of this paper, we have:

$$n_2^{1+\gamma}(\hat{d}_{2,av} - d^0) \xrightarrow{d} \frac{d_l + d_u}{2}.$$

It is also not difficult to see that this limit distribution is symmetric about 0.

Henceforth, the notation Argmin will denote the simple average of the minimal and maximal minimizers of a compound Poisson process. The quantity $|\alpha_0 - \beta_0|/\sigma$ will be denoted as SNR (signal-to-noise ratio). The higher the SNR, the more advantageous the estimation of the change-point at any given sample size would be. By (3), we have

$$\mathbb{M} \equiv \mathbb{M}_{|\alpha_0 - \beta_0|,\epsilon,C(K_1,\lambda_1,\gamma_1)} \stackrel{d}{\equiv} \frac{1}{C(K_1,\lambda_1,\gamma_1)} \operatorname{Arg\,min} \mathbb{M}_{\text{SNR},Z,1},$$

where $Z$ is the standard normal random variable. This is a consequence of the fact that $\epsilon/\sigma \sim N(0,1)$. From Theorem 2 and the above display, we have

$$n^{1+\gamma_1}(\hat{d}_{2,av} - d^0) \xrightarrow{d} \frac{1}{\lambda_2^{1+\gamma_1}} 2K_1 \left(\frac{\lambda_2}{\lambda_1}\right)^{\gamma_1} \operatorname{Arg\,min} \mathbb{M}_{\text{SNR},Z,1}$$

$$\stackrel{d}{\equiv} \frac{2K_1}{\lambda_2 \lambda_1^{\gamma_1}} \operatorname{Arg\,min} \mathbb{M}_{\text{SNR},Z,1}.$$

For a given sample size $n$, the objective is to choose the parameters $K_1$ and $\gamma_1$, so that the change-point $d^0$ would be contained in the prescribed interval $[\hat{d}_{1,av} - K_1 n_1^{-\gamma_1}, \hat{d}_{1,av} + K_1 n_1^{-\gamma_1}]$ with very high probability. In practice, this translates to choosing $K_1$ in such a way that $K_1 n_1^{-\gamma_1} \approx C_{\zeta_1}/n_1$, where $C_{\zeta_1}$ is the upper $\zeta_1$th quantile of the distribution of $\operatorname{Arg\,min} \mathbb{M}_{\text{SNR},Z,1}$, which is symmetric about 0. The value of $\zeta_1$ should be small, say 0.0005. In other words, our goal is to "zoom in," but not excessively so that we systematically start missing $d^0$ in the prescribed sampling interval. With this choice for $K_1$ we then have

(5)
$$\hat{d}_{2,av} - d^0 \approx \frac{2C_{\zeta_1}}{n_1^{1-\gamma_1} \lambda_2 \lambda_1^{\gamma_1} n^{-(1+\gamma_1)}} \operatorname{Arg\,min} \mathbb{M}_{\text{SNR},Z,1}$$
$$= \frac{2C_{\zeta_1}}{n^2 \lambda_2 \lambda_1} \operatorname{Arg\,min} \mathbb{M}_{\text{SNR},Z,1}.$$

It can be seen that the right-hand side is minimized when setting $\lambda_1 = \lambda_2 = 1/2$ (equal allocation of samples in the first and second stages).

Consider now a one-stage procedure with the covariates sampled from a density $f_X$, with the estimate of the change-point once again chosen to be the simple average of the minimal and maximal minimizers; call this $\hat{d}_{av}$.

In this case, the standard change-point asymptotics in conjunction with (2) and (3) give

$$n(\hat{d}_{av} - d^0) \xrightarrow{d} \frac{1}{f_X(d^0)} \operatorname{Arg\,min} \mathbb{M}_{\text{SNR},Z,1}.$$

This immediately provides an expression of the asymptotic efficiency of the two-stage procedure with respect to the one-stage (in terms of ratios of approximate standard deviations) given by

(6) $$\operatorname{ARE}_{2,1}(n) \approx \frac{n}{8 C_{\zeta_1} f_X(d^0)}.$$

It is not difficult to see that the same approximate formula for the ARE holds for some other measures of dispersion, besides the standard deviation. Let

$$\operatorname{ARE}_{2,1,\text{MAD}}(n) \equiv \frac{E|\hat{d}_{av} - d^0|}{E|\hat{d}_{2,av} - d^0|} \quad \text{and} \quad \operatorname{ARE}_{2,1,\text{IQR}}(n) \equiv \frac{\operatorname{IQR}(\hat{d}_{av})}{\operatorname{IQR}(\hat{d}_{2,av})},$$

where both first and second-stage estimates are based on samples of size $n$, and $\operatorname{IQR}(X)$ denotes the interquartile range of the distribution of a random variable $X$. Then, following similar steps to those involved in calculating the ARE based on standard deviations, we conclude that

$$\operatorname{ARE}_{2,1,\text{MAD}}(n) \approx \operatorname{ARE}_{2,1,\text{IQR}}(n) \approx \frac{n}{8 C_{\zeta_1} f_X(d^0)}.$$

The accuracy of the above approximation is confirmed empirically through a simulation study. The setting involves a change-point model given by

(7) $$y_i = 0.5 I(x_i \leq 0.5) + 1.5 I(x_i > 0.5) + \epsilon_i, x_i \in (0,1).$$

The variance $\sigma^2$ was chosen so that the SNR defined as $(\beta_0 - \alpha_0)/\sigma = 1, 2, 5$ and 8 and the sample size varies in increments of 50 from 50 to 1500. The results based on an interval corresponding to $\zeta_1 = 0.0025$ and 5000 replications are shown in Figure 2. Further, for the two-stage procedure the smallest and largest 5 "outliers" were dropped, since they introduced too much variability, especially for the standard deviation based ARE. These "outliers" correspond to the cases where the true parameters $d^0$ was not contained in the zoom-in sampling interval. It can be seen that there is great agreement between the theoretical formula for the ARE and the empirical ARE, for all three performance measures employed.

REMARK. The formula for the ARE in (6) says that the "agnostic" two-stage procedure ("agnostic" since the covariates are sampled uniformly at each stage) will eventually, that is, with increasing $n$, surpass any one stage procedure, no matter the amount of background information incorporated about the location of the change-point in the one-stage process. One can think of an "oracle-type" one-stage procedure where the experimenter samples the covariates from a density that peaks in a neighborhood of $d^0$ relative to the uniform density [corresponding to high values of $f_X(d^0)$]. The faster convergence rate of the two-stage procedure relative to this one-stage procedure guarantees that with increasing $n$, the ARE will always go to infinity. Further, expression (6) provides an approximation to the minimal sample size required for the two-stage procedure to outperform the "classical" one, a result verified through simulations (not shown here).

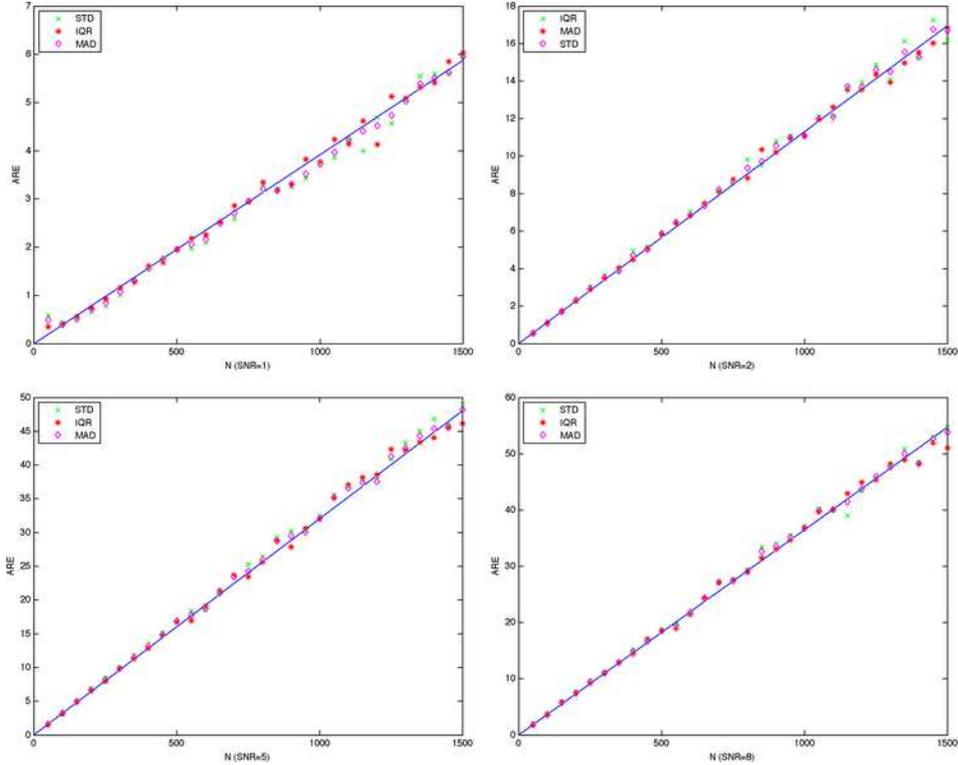

Fig. 2. *Top panels: ARE for standard deviation, IQR and MAD measures for* SNR = 1 *(left) and 2 (right). Bottom panels: corresponding ARE for* SNR = 5 *(left) and 8 (right).*

REMARK. A uniform density has been considered up to this point for sampling second-stage, covariate-response pairs. We examine next the case of using an arbitrary sampling density $g_U(\cdot)$ supported on the interval $[\hat{d}_{n_1} - Kn_1^{-\gamma}, \hat{d}_{n_1} + Kn_1^{-\gamma}]$ and symmetric around $\hat{d}_{n-1}$. A natural choice for such a density is $g_U(d^0) = h((d^0 - \hat{d}_{n_1})/n_1^{-\gamma}K)(n_1^\gamma/K)$, for a density $h(\cdot)$ supported on $[-1, 1]$ and symmetric about 0. Analogous arguments to those used in the proof of Theorem 2 establish that the random vector $n_2^{1+\gamma}(\hat{d}_{n_2,l} - d^0, \hat{d}_{n_2,u} - d^0)$ converges in distribution to

$$(d_l(|\alpha_0 - \beta_0|, \epsilon, C(K, \lambda, \gamma, h)), d_u(|\alpha_0 - \beta_0|, \epsilon, C(K, \lambda, \gamma, h))),$$

where $C(K, \lambda, \gamma, h) = (\lambda/(1-\lambda))^\gamma (h(0)/K)$. With the error term normally distributed, as assumed in this section, we get

$$n^{1+\gamma}(\hat{d}_{2,av} - d^0) \xrightarrow{d} \frac{K}{h(0)(1-\lambda)\lambda^\gamma} \operatorname{Arg\,min} \mathbb{M}_{\mathrm{SNR},Z,1}$$

and it can be readily checked that the approximate ARE formula reduces to $\mathrm{ARE}_{2,1}(n) \approx nh(0)/(4C_\zeta f_X(d^0))$. It can be seen that the more "peaked" the sampling density $g_U$ (equivalently $h$) the greater the efficiency gains. However, one needs to be careful, since the above formula is obtained through asymptotic considerations. In finite samples, a very peaked density around $\hat{d}_{n_1}$ may not perform well, since bias issues [involving $(d^0 - \hat{d}_{n_1})$] must also be taken into account.

REMARK. Simulation results indicate that in the presence of a small budget of available points ($n = 20$ or $50$), the efficiency of the two-stage estimator can be improved by employing a uniform (equispaced) design in the first-stage. The reason is that such a design reduces the sampling variability of the covariate $x$, which leads to improved localization of the change-point. However, the approximate formula for the ARE discussed above is no longer valid when we use a uniform design at the first-stage, since the asymptotics of the one-stage estimator are then no longer described by the minimizer of a compound Poisson process. For further discussion on uniform sampling designs, see Section 4.1 of Lan, Banerjee and Michailidis (2007).

We turn our attention to the allocation of parameters $K_i, \lambda_i$ and $\gamma_i, 1 \leq i \leq P-1$ for the $P$ stage procedure. We start with a three-stage procedure and generalize afterward. From the theoretical result of Section 4.1 we have

$$n_3^{2+\gamma_2}(\hat{d}_{n_3,av} - d^0) \xrightarrow{d} 2K_2\left(\frac{\lambda_3}{\lambda_2}\right)^{1+\gamma_2} \operatorname{Arg\,min} \mathbb{M}_{\mathrm{SNR},Z,1} \quad (8)$$

$$\implies n^{2+\gamma_2}(\hat{d}_{n_3,av} - d^0)$$

$$\xrightarrow{d} \frac{2K_2}{\lambda_3 \lambda_2^{1+\gamma_2}} \operatorname{Arg\,min} \mathbb{M}_{\mathrm{SNR},Z,1} \quad (9)$$

or

$$\hat{d}_{n_3,av} - d^0 \stackrel{d}{\approx} \frac{2K_2}{n^{2+\gamma_2} \lambda_3 \lambda_2^{1+\gamma_2}} \operatorname{Arg\,min} \mathbb{M}_{\mathrm{SNR},Z,1}.$$

Once again, our objective is to choose $K_2$ and $\gamma_2$ so that for any fixed sample size $n$, with very high probability $d^0$ would be contained in the sampling interval $[\hat{d}_{n_2,av} - K_2 n_2^{-(1+\gamma_2)}, \hat{d}_{n_2,av} + K_2 n_2^{-(1+\gamma_2)}]$. Thus, using (5) we require that $K_2 n_2^{-(1+\gamma_2)} \approx (2C_{\zeta_1} C_{\zeta_2})/(n^2 \lambda_2 \lambda_1)$, or $K_2 = (2C_{\zeta_1} C_{\zeta_2})/(\lambda_1 \lambda_2^{-\gamma_2} n^{1-\gamma_2})$, with both $\zeta_1, \zeta_2$ being very small. With this choice for $K_2$ we obtain

$$\hat{d}_{n_3,av} - d^0 \stackrel{d}{\approx} \frac{4 C_{\zeta_1} C_{\zeta_2}}{n^3 \lambda_1 \lambda_2 \lambda_3} \operatorname{Arg\,min} \mathbb{M}_{\mathrm{SNR},Z,1}.$$

It is then easy to see that the right-hand size is minimized when $\lambda_1 = \lambda_2 = \lambda_3 = 1/3$. The ARE formula for the three-stage procedure under the proposed (equal) allocation of samples was validated through a simulation based on 5000 replications. The same stump model was employed and the values for SNR were set equal to 5 and 8, while those for $\zeta_1 = \zeta_2 = 0.0025$. The results using a similar trimming of "outliers" for the estimates of the three-stage procedure (3 smallest and 3 largest in this case) are shown in Figure 3 and exhibit great agreement with the theoretical formula.

Proceeding analogously to the three-stage procedure, we obtain that at stage $q-1$,

$$\hat{d}_{n_{q-1},av} - d^0 \stackrel{d}{\approx} \frac{2^{q-2} C_{\zeta_1} \cdots C_{\zeta_{q-2}}}{n^{q-1} \lambda_1 \cdots \lambda_{q-1}} \operatorname{Arg\,min} \mathbb{M}_{\mathrm{SNR},Z,1}.$$

A level $1 - 2\zeta_{q-1}$ C.I. for $d^0$ (for a small $\zeta_{q-1}$) is then given by

$$(10) \quad \left[ \hat{d}_{n_{q-1},av} - \frac{2^{q-2} C_{\zeta_1} \cdots C_{\zeta_{q-2}} C_{\zeta_{q-1}}}{n^{q-1} \lambda_1 \cdots \lambda_{q-1}}, \hat{d}_{n_{q-1},av} + \frac{2^{q-2} C_{\zeta_1} \cdots C_{\zeta_{q-2}} C_{\zeta_{q-1}}}{n^{q-1} \lambda_1 \cdots \lambda_{q-1}} \right].$$

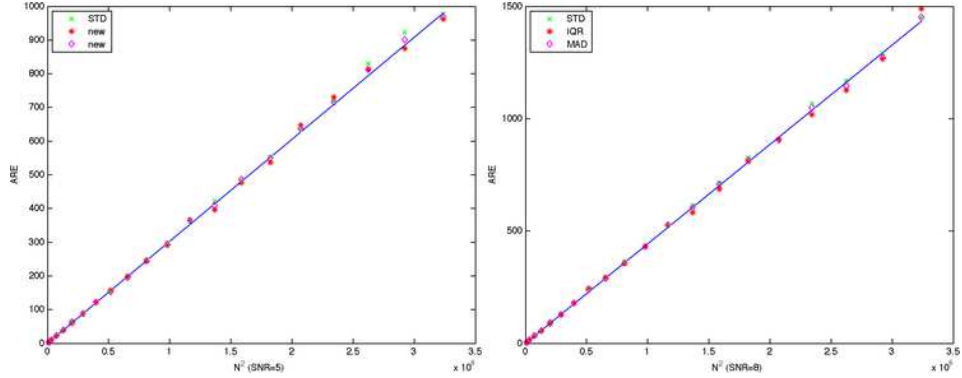

Fig. 3. *ARE for standard deviation, IQR and MAD measures for* $\mathrm{SNR} = 5$ *(left) and 8 (right).*

The sampling neighborhood of $\hat{d}_{n_{q-1},av}$ from which we sample at stage $q$ is formally given by

$$\tag{11} \left[\hat{d}_{n_{q-1},av} - \frac{K_{q-1}}{n_{q-1}^{q-2+\gamma_{q-1}}}, \hat{d}_{n_{q-1},av} + \frac{K_{q-1}}{n_{q-1}^{q-2+\gamma_{q-1}}}\right]$$

and yields that

$$(\hat{d}_{n_q,av} - d^0) \stackrel{d}{\approx} \frac{2K_{q-1}}{\lambda_q \lambda_{q-1}^{(q-2)+\gamma_{q-1}} n^{(q-1)+\gamma_{q-1}}} \operatorname{Arg\,min} \mathbb{M}_{SNR,Z,1}.$$

Equating the neighborhoods (10) and (11) gives

$$K_{q-1} = \frac{2^{q-2} C_{\zeta_1} \cdots C_{\zeta_{q-2}} C_{\zeta_{q-1}}}{n^{1-\gamma_{q-1}} \lambda_1 \cdots \lambda_{q-2} \lambda_{q-1}^{-\gamma_{q-1}-q+3}}.$$

Plugging this choice of $K_{q-1}$ into the expression for the approximate distribution of $\hat{d}_{n_q,av} - d^0$ gives

$$\hat{d}_{n_q,av} - d^0 \stackrel{d}{\approx} \frac{2^{q-1} C_{\zeta_1} \cdots C_{\zeta_{q-2}} C_{\zeta_{q-1}}}{n^q \lambda_1 \cdots \lambda_{q-1} \lambda_q} \operatorname{Arg\,min} \mathbb{M}_{SNR,Z,1}.$$

Therefore, the $P$ stage estimate $(\hat{d}_{n_p} - d^0)$ is approximately distributed as

$$\frac{2^{P-1} C_{\zeta_1} \cdots C_{\zeta_{P-1}}}{\lambda_1 \cdots \lambda_P n^P} \operatorname{Arg\,min} \mathbb{M}_{SNR,Z,1},$$

which shows that an equal allocation of samples is warranted (i.e., $\lambda_i = 1/P, 1 \leq i \leq P$). Some straightforward algebra establishes that

$$\operatorname{ARE}_{P,1}(n) \approx \frac{n^{P-1} \lambda_1 \lambda_2 \cdots \lambda_P}{2^{P-1} C_{\zeta_1} \cdots C_{\zeta_{P-1}} f_X(d^0)} = \frac{n^{P-1}}{2^{P-1} P^P f_X(d^0) C_{\zeta_1} \cdots C_{\zeta_{P-1}}},$$

where $f_X(d^0)$ has the same connotation as in (6).

REMARK. An interesting question that arises in practice is that of how many stages to use given a fixed budget of $n$ points. Obviously the answer depends on the underlying SNR, but $n/P$ should be fairly large. Extensive numerical work with two, three and four-stage procedures indicates that $n/P$ in the range of 35–50 points performs well. A related issue is how to choose the values of $\zeta_i$ that determine the coverage of the successive neighborhoods. Notice that with probability $1-(1-2\zeta_1)$, $d^0$ is *not* contained in the sampling interval neighborhood around $\hat{d}_{n_1}$. Proceeding analogously and due to the independence of sampling at subsequent stages, it can be seen that with probability $1 - \prod_{i=1}^{P-1}(1 - 2\zeta_i)$ the shrinking sampling intervals are going to miss $d^0$ at some stage. Taking all the $2\zeta_i$'s to be equal to $\psi$, we get that this probability is given by $1 - (1-\psi)^{P-1}$. Setting the probability of trapping $d^0$ after $P$ stages equal to $(1-\delta)$ ($\delta$ very small) we get that $\psi = 1 - (1-\delta)^{1/(P-1)}$, which provides a good guideline for determining the size of the sampling neighborhoods. There remains the issue of the sampling interval not containing $d^0$ at some stage. The ARE results based on simulations for a 3-stage procedure, that are shown in Figure 3, are somewhat optimistic due to the trimming of outliers discussed previously; however, even with outliers included, improvements still occur, albeit small compared to the one-stage procedure.

**5. Confidence intervals for the change-point.** We compare next the performance of exact confidence intervals based on the result established in Theorem 2 to those proposed in Ferger (2004). Moreover, confidence intervals for finite samples will be constructed following the discussion in Section 4.

For all these comparisons, simulations were run for a stump model with $\alpha_0 = 0.5, \beta_0 = 1.5, d^0 = 0.5$ and sample sizes $n = 50, 100, 200, 500, 1000$ with $N = 2000$ replicates for each $n$. Confidence intervals for $d^0$ based on the minimal minimizer $\hat{d}_{n_2,l}$, the maximal minimizer $\hat{d}_{n_2,u}$ and the average minimizer $\hat{d}_{n_2,av} = (\hat{d}_{n_2,l} + \hat{d}_{n_2,u})/2$ were constructed. Two values of $\gamma = 1/2$ and $2/3$ and two values of $K = 1$ and $2$ were used together with the optimal allocation $\lambda \equiv \gamma/(1+\gamma)$ as discussed in Section 4. The confidence level was set at $1 - \tau = 0.95$ and the percentage of replicates for which the true change-point was included in the corresponding intervals, as well as the average length of each interval, were recorded. In what follows, the symbols $d_l$ and $d_u$ have the same connotations as in the proof of Theorem 3.2.

5.1. *Conservative intervals.* Using the results of Ferger (2004), based on any two-stage estimator $\hat{d}_{n_2}$, we propose an asymptotically conservative confidence interval for $d^0$ at level $1 - \tau$:

$$I_{n_2}(\tau) := (\hat{d}_{n_2} - b/n_2^{1+\gamma}, \hat{d}_{n_2} - a/n_2^{1+\gamma}),$$

where $a < b$ are any solutions of the inequality

$$\text{Prob}(d_u < b) - \text{Prob}(d_l \leq a) \geq 1 - \tau.$$

Based on the smallest, largest and average minimizers at stage two, we therefore obtain intervals

$$I_{n_2,l}(\tau) = (\hat{d}_{n_2,l} - b/n_2^{1+\gamma}, \hat{d}_{n_2,l} - a/n_2^{1+\gamma}),$$
$$I_{n_2,u}(\tau) = (\hat{d}_{n_2,u} - b/n_2^{1+\gamma}, \hat{d}_{n_2,u} - a/n_2^{1+\gamma})$$

and

$$I_{n_2,av}(\tau) = (\hat{d}_{n_2,av} - b/n_2^{1+\gamma}, \hat{d}_{n_2,av} - a/n_2^{1+\gamma}),$$

where $a$ is the $\tau/2$th quantile of $d_l$ and $b$ is the $(1-\tau/2)$th quantile of $d_u$. At this point, these quantiles do not seem to be analytically determinable but can certainly be simulated to a reasonable degree of approximation.

In Table 1 the coverage probabilities together with the length of the confidence intervals are shown for a number of combinations of sample sizes and tuning parameters and with the SNR set equal to 5. It can be seen that the recorded coverage exceeds the nominal level of 95% and almost approaching perfect (100%) coverage for the average minimizer.

5.1.1. *Exact confidence intervals.* On the other hand, since Theorem 3.2 provides us with the *exact* asymptotic distributions of the sample minimizers, we can construct asymptotically exact (level $1-\tau$ confidence intervals) as follows:

$$\tilde{I}_{n_2,l}(\tau) = (\hat{d}_{n_2,l} - b_l/n_2^{1+\gamma}, \hat{d}_{n_2,l} - a_l/n_2^{1+\gamma}),$$
$$\tilde{I}_{n_2,u}(\tau) = (\hat{d}_{n_2,u} - b_u/n_2^{1+\gamma}, \hat{d}_{n_2,u} - a_u/n_2^{1+\gamma}),$$
$$\tilde{I}_{n_2,av}(\tau) = (\hat{d}_{n_2,av} - b_{av}/n_2^{1+\gamma}, \hat{d}_{n_2,av} - a_{av}/n_2^{1+\gamma}),$$

where $a_l, b_l, a_u, b_u, a_{av}$ and $b_{av}$ are the exact quantiles [$a_l, a_u$ and $a_{av}$ correspond to $\tau/2$th quantiles and $b_l, b_u$ and $b_{av}$ correspond to $(1-\tau/2)$th quantiles] of $d_l$, $d_u$ and $(d_l + d_u)/2$, respectively.

In Table 2 the coverage probabilities together with the length of the confidence intervals are shown for a number of combinations of sample sizes and tuning parameters and with the SNR set equal to 5. It can be seen that the coverage probabilities are fairly close to their nominal values, especially for $\gamma = 2/3$. Further, their length is almost half of those obtained according to Ferger's (2004) method. Finally, it should be noted that analogous results were obtained for $\text{SNR} = 2$ and 8 (not shown due to space considerations).

5.1.2. *Construction of confidence intervals in finite samples.* Confidence intervals in finite samples can also be based on the adaptive parameter allocation strategies discussed in Section 4 for the two-stage procedure. We briefly discuss this below, adopting notation from that section.

From (5) with equal allocation of points between the two stages we have that

$$\hat{d}_{n_2,av} - d^0 \stackrel{d}{\approx} \frac{8C_{\zeta_1}}{n^2} \operatorname{Arg\,min} \mathbb{M}_{\text{SNR},Z,1}.$$

Therefore, an approximate level $1-\tau$ confidence interval is given by

(12) $$\left[\widehat{d_{n_2,av}} - \frac{8C_{\zeta_1}C_{\tau/2}}{n^2}, \widehat{d_{n_2,av}} + \frac{8C_{\zeta_1}C_{\tau/2}}{n^2}\right].$$

Simulations were run for the above stump model for four different sample sizes: 50, 100, 200 and 500 with 5000 replicates for each sample size and for three different values of SNR = 2, 5, 8. Confidence intervals as defined above were constructed (with $\tau = 0.05$). The percentage of intervals containing the true change-point together with their length were recorded and shown in Table 3.

TABLE 1
*95% conservative confidence intervals for a combination of sample sizes and the tuning parameters $\gamma, K$ and for $\mathrm{SNR} = 5$*

|  | $n = 50$ | | $n = 100$ | | $n = 200$ | | $n = 500$ | | $n = 1000$ | |
| --- | --- | --- | --- | --- | --- | --- | --- | --- | --- | --- |
|  | $K=1$ | $K=2$ | $K=1$ | $K=2$ | $K=1$ | $K=2$ | $K=1$ | $K=2$ | $K=1$ | $K=2$ |
|  | | | | | $\gamma = \frac{1}{2}$ | | | | | |
| $\hat{I}_{n_2,l}$ | 97.40% | 97.55% | 97.40% | 98.20% | 97.65% | 98.00% | 96.55% | 97.05% | 97.15% | 97.75% |
|  | (0.0580) | (0.1208) | (0.0205) | (0.0427) | (0.0072) | (0.0151) | (0.0018) | (0.0038) | (0.0006) | (0.0014) |
| $\hat{I}_{n_2,u}$ | 97.05% | 98.60% | 97.65% | 97.05% | 97.65% | 97.85% | 97.40% | 97.90% | 97.80% | 98.00% |
|  | (0.0580) | (0.1208) | (0.0205) | (0.0427) | (0.0072) | (0.0151) | (0.0018) | (0.0038) | (0.0006) | (0.0014) |
| $\hat{I}_{n_2,av}$ | 99.80% | 99.95% | 99.80% | 99.95% | 100% | 100% | 99.80% | 100% | 99.70% | 99.95% |
|  | (0.0580) | (0.1208) | (0.0205) | (0.0427) | (0.0072) | (0.0151) | (0.0018) | (0.0038) | (0.0006) | (0.0014) |
|  | | | | | $\gamma = \frac{2}{3}$ | | | | | |
| $\hat{I}_{n_2,l}$ | 98.15% | 97.70% | 97.65% | 98.30% | 97.90% | 97.70% | 97.65% | 97.30% | 97.75% | 97.70% |
|  | (0.0299) | (0.0581) | (0.0094) | (0.0183) | (0.0030) | (0.0058) | (0.0006) | (0.0013) | (0.0002) | (0.0004) |
| $\hat{I}_{n_2,u}$ | 98.00% | 98.20% | 97.85% | 98.55% | 97.90% | 98.10% | 98.30% | 97.75% | 97.60% | 98.50% |
|  | (0.0299) | (0.0581) | (0.0094) | (0.0183) | (0.0030) | (0.0058) | (0.0006) | (0.0013) | (0.0002) | (0.0004) |
| $\hat{I}_{n_2,av}$ | 99.60% | 99.95% | 99.60% | 99.90% | 99.85% | 99.95% | 99.85% | 99.90% | 99.85% | 99.95% |
|  | (0.0299) | (0.0581) | (0.0094) | (0.0183) | (0.0030) | (0.0058) | (0.0006) | (0.0013) | (0.0002) | (0.0004) |

TABLE 2
*95% exact confidence intervals for a combination of sample sizes and tuning parameters $\gamma, K$ for SNR $= 5$*

|  | $n = 50$ | | $n = 100$ | | $n = 200$ | | $n = 500$ | | $n = 1000$ | |
| --- | --- | --- | --- | --- | --- | --- | --- | --- | --- | --- |
|  | $K = 1$ | $K = 2$ | $K = 1$ | $K = 2$ | $K = 1$ | $K = 2$ | $K = 1$ | $K = 2$ | $K = 1$ | $K = 2$ |
|  | | | | | $\gamma = \frac{1}{2}$ | | | | | |
| $\tilde{I}_{n_2,l}$ | 95.00% | 94.20% | 93.80% | 95.25% | 95.25% | 94.10% | 93.40% | 94.50% | 95.05% | 94.90% |
|  | (0.0283) | (0.0599) | (0.0100) | (0.0212) | (0.0035) | (0.0075) | (0.0009) | (0.0019) | (0.0003) | (0.0007) |
| $\tilde{I}_{n_2,u}$ | 94.20% | 96.50% | 94.80% | 94.95% | 94.45% | 95.85% | 94.85% | 95.90% | 95.50% | 95.85% |
|  | (0.0294) | (0.0602) | (0.0104) | (0.0213) | (0.0037) | (0.0075) | (0.0009) | (0.0019) | (0.0003) | (0.0007) |
| $\tilde{I}_{n_2,av}$ | 94.30% | 95.25% | 94.35% | 94.05% | 94.40% | 95.45% | 93.20% | 94.85% | 94.55% | 95.30% |
|  | (0.0236) | (0.0487) | (0.0083) | (0.0172) | (0.0029) | (0.0061) | (0.0007) | (0.0015) | (0.0003) | (0.0005) |
|  | | | | | $\gamma = \frac{2}{3}$ | | | | | |
| $\tilde{I}_{n_2,l}$ | 95.65% | 95.40% | 95.05% | 96.05% | 95.35% | 95.45% | 95.30% | 95.30% | 95.05% | 95.90% |
|  | (0.0148) | (0.0277) | (0.0047) | (0.0087) | (0.0015) | (0.0027) | (0.0003) | (0.0006) | (0.0001) | (0.0002) |
| $\tilde{I}_{n_2,u}$ | 95.40% | 96.00% | 95.65% | 96.60% | 95.80% | 96.15% | 96.20% | 96.60% | 95.15% | 96.85% |
|  | (0.0149) | (0.0302) | (0.0047) | (0.0095) | (0.0015) | (0.0030) | (0.0003) | (0.0006) | (0.0001) | (0.0002) |
| $\tilde{I}_{n_2,av}$ | 95.15% | 96.45% | 95.20% | 96.95% | 94.75% | 96.10% | 95.30% | 96.15% | 94.05% | 96.10% |
|  | (0.0120) | (0.0253) | (0.0038) | (0.0080) | (0.0012) | (0.0025) | (0.0003) | (0.0005) | (0.0001) | (0.0002) |

TABLE 3
*95% confidence intervals constructed using the adaptive parameter allocation strategy for different sample sizes and SNR with $\zeta_1 = 0.0005$*

|   | SNR = 2 | | SNR = 5 | | SNR = 8 | |
| --- | --- | --- | --- | --- | --- | --- |
| $N$ | Coverage | Length | Coverage | Length | Coverage | Length |
| 50  | 93.24% | 0.2780 | 95.48% | 0.0383  | 95.68% | 0.0329  |
| 100 | 94.08% | 0.0695 | 95.24% | 0.0096  | 95.54% | 0.0082  |
| 200 | 94.48% | 0.0174 | 94.78% | 0.0024  | 95.16% | 0.0021  |
| 500 | 94.82% | 0.0028 | 95.08% | 0.00038 | 94.94% | 0.00033 |

We examine next the performance of confidence intervals in finite samples, but where a uniform (equispaced) design is used in the first-stage (results shown in Table 4) and in both stages (results shown in Table 5). The setting is identical to that used in Table 3. It is not clear how the tuning parameter $C_{\zeta_1}$, that determines the interval from which sampling is done at the second stage, should be chosen in this case, since the first-stage estimate may not even have an asymptotic distribution in the proper sense. Therefore, the same $C_{\zeta_1}$ value as the one used in Table 3 was employed. It can be seen that a uniform design used in the first-stage does not improve performance in terms of coverage or length. However, using a uniform design in both stages and setting $C_{\zeta_1}$ and $C_{\tau/2}$ to the same values as in Table 3 leads to rather conservative confidence intervals, especially for larger sample sizes and higher values of SNR. Notice that the lengths of the confidence intervals are identical to those in Table 3 due to the choice of the tuning parameters $C_{\zeta_1}$ and $C_{\tau/2}$. Nevertheless, experience shows that a uniform design used in the first-stage gives better mean squared errors in small samples, or when $d^0$ is closer to the boundary of the covariate's support.

TABLE 4
*95% confidence interval constructed using the adaptive parameter allocation strategy for different sample sizes and SNR with $\zeta_1 = 0.0005$ using a uniform design in the first-stage*

|   | SNR = 2 | | SNR = 5 | | SNR = 8 | |
| --- | --- | --- | --- | --- | --- | --- |
| $N$ | Coverage | Length | Coverage | Length | Coverage | Length |
| 50  | 93.72% | 0.2780 | 95.14% | 0.0383  | 95.56% | 0.0329  |
| 100 | 93.88% | 0.0695 | 95.12% | 0.0096  | 95.20% | 0.0082  |
| 200 | 94.62% | 0.0174 | 95.52% | 0.0024  | 95.52% | 0.0021  |
| 500 | 94.72% | 0.0028 | 94.96% | 0.00038 | 95.12% | 0.00033 |

TABLE 5
*95% confidence interval constructed using the adaptive parameter allocation strategy for different sample sizes and SNR with $\zeta_1 = 0.0005$ using a uniform design in both stages*

| $N$ | SNR = 2 | | SNR = 5 | | SNR = 8 | |
|---|---|---|---|---|---|---|
| | Coverage | Length | Coverage | Length | Coverage | Length |
| 50  | 95.06% | 0.2780 | 100.00% | 0.0383  | 100.00% | 0.0329  |
| 100 | 96.66% | 0.0695 | 99.98%  | 0.0096  | 100.00% | 0.0082  |
| 200 | 96.94% | 0.0174 | 100.00% | 0.0024  | 100.00% | 0.0021  |
| 500 | 97.32% | 0.0028 | 99.96%  | 0.00038 | 100.00% | 0.00033 |

5.2. *Data application.* We revisit the motivating application and estimate the change-point using both the "classical" and the developed two-stage procedures. The total budget was set to $n = 70$ and the model fitted to the natural logarithm of the delays comprised two linear segments with a discontinuity. Given that the data (134 system loadings and their corresponding average delays) have been collected in advance, a sampling mechanism close in spirit to selecting covariate values from a uniform distribution was employed for both procedures. Specifically, the necessary number of points was drawn from a uniform distribution in the $[0, 1]$ interval and amongst the available 134 loadings the ones closest to the sampled points were selected, together with their corresponding responses. An analogous strategy was used when a uniform design was employed in the first-stage of the adaptive procedure. For the two-stage procedure, we set $\lambda = 1/2$ and the remaining tuning parameters to those values provided by the adaptive strategy discussed in Section 4, with $\zeta_1 = 0.0005$. The results of the "classical" procedure, the two-stage adaptive procedure with sampling from a uniform distribution in both stages and from a uniform design in the first-stage and the uniform distribution in the second stage are depicted in the left, center and right panels of Figure 4, respectively. The depicted fitted regression models are based on the first-stage estimates for the two-stage procedure. Further, the sampled points from the two stages are shown as solid (first-stage) and open (second stage) circles. It can be seen that the heavier sampling in the neighborhood of the first-stage estimate of the change-point improves the estimate given the available evidence from all 134 points shown in Figure 1.

The estimated change-point from the "classical" procedure is $\hat{d}_n = 0.737$ with a 95% confidence interval $(0.682, 0.793)$. Using a uniform distribution in both stages gave an estimate $\hat{d}_{n_2} = 0.796$ with a 95% confidence interval $(0.781, 0.811)$. On the other hand, a combination of a uniform design in the first-stage with that of a uniform distribution in the second stage yielded an estimate $\hat{d}_{n_2} = 0.802$ with a 95% confidence interval $(0.787, 0.817)$. As shown in this case and validated through other data examples, the use of uniform design in the first-stage proves advantageous in practice, especially for small samples or in situations where the discontinuity lies fairly close to the boundary of the design region.

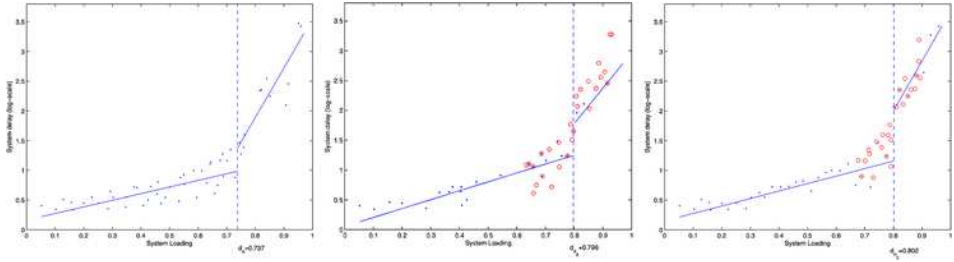

Fig. 4. *Sampled points (from first-stage solid circles and from second-stage open circles) together with the fitted parametric models and estimated change-point, based on a total budget of $n = 70$ points, obtained from the "classical" procedure (left panel), the two-stage adaptive procedure with sampling from a uniform distribution in both stages (center panel) and from a uniform design in the first-stage and the uniform distribution in the second stage (right panel).*

**6. Concluding remarks.** In this study, a multistage adaptive procedure for estimating a jump discontinuity in a parametric regression model was introduced and its properties investigated. Specifically, it was established that the rate of convergence of the least squares estimator of the change-point can be accelerated and its asymptotic distribution derived. Several issues pertaining to the tuning of the parameters involved in the procedure were examined and resolved. In practice, it is generally recommended that in the presence of a small budget of points, a uniform design in the first-stage be employed. At present the parameters of the parametric model are estimated in the first-stage of the experiment and held fixed in subsequent stages. One may wonder why the parameters are not re-estimated in the presence of additional samples. The main reason is that the additional samples obtained from a shrinking neighborhood around $d^0$ are not particularly helpful for estimating global regression parameters. The sole exception is the stump model, where using all $n$ points provides better estimates, especially for small budgets.

The results have been established under the assumption of a bounded covariate, because this is usually the case in applications. However, the accelerated convergence can be achieved with an unbounded covariate as well, so long as the first-stage estimate is consistent at rate $n$, which will be the case under fairly mild conditions.

Another issue is the case where there are multiple change-points, but with their number known a priori. In such a case, one will estimate each of the change-points at stage one (at rate $n$) and subsequently construct appropriate neighborhoods around these estimates. Notice that asymptotically the neighborhoods are disjoint and therefore one can perform the second-stage experiment for a single change-point model *independently* on each of them. However, there is an issue of allocation of second-stage points in the neighborhoods. Some algebra involving minimization of the sum of the asymptotic variances of the second stage estimators shows that 50% of the budget should be allocated at the first-stage, while the remaining budget should be allocated proportionally to the limiting standard deviations of the first-stage estimates, which can be estimated from the first-stage experiment. The real challenge with multiple change-points is when their number is unknown and needs to be determined from the data. This topic remains to be studied.

We briefly address the case of multiple covariates next. The simplest possible model is for two covariates $x^1, x^2$ on the unit square and defined as follows:

$$Y_i = \alpha_1 1(x_i^1 \leq d^1) + \beta_1 1(x_i^1 > d^1) + \alpha_2 1(x_i^2 \leq d^2) + \beta_2 1(x_i^2 > d^2) + \epsilon_i.$$

For the sake of simplicity, suppose that a budget of $2n$ points is available. The change-points $d^1$ and $d^2$ can be estimated at rate $n$, by sampling $n$ points from a uniform distribution on the unit square and solving the corresponding least squares problem. One can then form neighborhoods around $\hat{d}^1$ and $\hat{d}^2$, respectively, of the type employed in the univariate models. Hence, a "shrunken" rectangular neighborhood of $(\hat{d}^1, \hat{d}^2)$ is obtained, from which one can sample $n$ points uniformly at the second stage. The least squares estimates from stage two will inherit the accelerated rate of convergence. Similar considerations apply to general parametric models.

Note that in this setting, there are changing regression models on each quadrant of $\mathbb{R}^2$ defined by two orthogonally intersecting lines. A more complex model arises when the lines intersect at an unknown angle. However, a treatment of such models is beyond the scope of this paper and is left as a topic of future research.

Finally, it should be noted that a multistage adaptive procedure would also work in the context of a nonparametric model with a jump discontinuity and produce analogous accelerations to the convergence rate of the employed estimator; a comprehensive treatment of this topic is currently under study.

## APPENDIX

PROOF OF LEMMA 3.2.  We first note that $\mathcal{D}_C^0$ (which we view as a metric subspace of $D[-C, C] \times D[-C, C]$) is a measurable subset of $D[-C, C] \times$ $D[-C,C]$. To establish convergence in distribution in the space $\mathcal{D}_C^0$, it therefore suffices to establish convergence in distribution in the larger space $D[-C,C] \times D[-C,C]$ [see the discussion in Example 3.1 of Billingsley (1999)]. This can be achieved by (a) establishing finite-dimensional convergence: showing that $\{\mathbb{M}_{n_2}(h_i), \mathbb{J}_{n_2}(h_i)\}_{i=1}^l \to \{\mathbb{M}(h_i), \mathbb{J}(h_i)\}_{i=1}^l$ for all $h_1, h_2, \ldots, h_l$ in $[-C, C]$ and (b) verifying tightness of the processes $(\mathbb{M}_{n_2}(h), \mathbb{J}_{n_2}(h))$ under the product topology. But this boils down to verifying marginal tightness.

Let $L^+(t) = E[e^{it(W-(\alpha_0+\beta_0)/2)}|U=d^{0+}] \equiv \lim_{d \to d^{0+}} E[e^{it(W-(\alpha_0+\beta_0)/2)}|U=d]$ and $L^-(t) = E[e^{it(W-(\alpha_0+\beta_0)/2)}|U=d^0]$. It is not difficult to see that $L^+$ is the characteristic function of the $V_i^+$'s while $L^-$ is the characteristic function of the $V_i^-$'s. In order to establish finite-dimensional convergence, we first show the convergence of one-dimensional marginals: that is, for a fixed $s$, $\mathbb{M}_{n_2}(s)$ converges in distribution to $\mathbb{M}(s)$. We do this via characteristic functions. Consider $\phi_s(t)$, the characteristic function of $\mathbb{M}(s)$ (with $s > 0$). We have

$$\begin{aligned}
E[e^{it\mathbb{M}_1(s)}] &= E[e^{it\sum_{0 \le k \le \nu^+(s)} V_k^+}] \\
&= \sum_{l=0}^{\infty} E[e^{it(V_1^+ + \cdots + V_l^+)}] \frac{e^{-(s/2K)\lambda/(1-\lambda)^\gamma}(s/(2K)\lambda/(1-\lambda)^\gamma)^l}{l!} \\
&= \sum_{l=0}^{\infty} \frac{(L^+(t))^l (s/(2K)\lambda/(1-\lambda)^\gamma)^l}{l!} e^{-(s/2K)\lambda/(1-\lambda)^\gamma} \\
&= e^{-(s/2K)\lambda/(1-\lambda)^\gamma} e^{L(t)(s/2K)\lambda/(1-\lambda)^\gamma} = e^{-(s/2K)\lambda/(1-\lambda)^\gamma(1-L(t))}.
\end{aligned}$$

We show that $Q_{n_2,s}(t) \equiv E[e^{it\mathbb{M}_{n_2}(s)}]$ converges to $\phi_s(t)$. Let $\xi_{n_1} = n_1(\hat{d}_{n_1} - d^0), \eta_{n_1,1} = \sqrt{n_1}(\hat{\alpha}_{n_1} - \alpha_0), \eta_{n_1,2} = \sqrt{n_1}(\hat{\beta}_{n_1} - \beta_0)$. We have

$$Q_{n_2,s}(t) = \int Q_{n_2,s}^\star(t, \eta_1, \eta_2, \xi) \, dZ_{n_1}(\eta_1, \eta_2, \xi),$$

where $Z_{n_1}$ is the joint distribution of $(\eta_{n_1,1}, \eta_{n_1,2}, \xi_{n_1})$ and

$$\begin{aligned}
Q_{n_2,s}^\star&(t, \eta_1, \eta_2, \xi) \\
&= E[e^{it\mathbb{M}_{n_2}^+(s)}|\eta_{n_1,1} = \eta_1, \eta_{n_1,2} = \eta_2, \xi_{n_1} = \xi] \\
&= E[e^{it(W_1 - (\hat{\alpha}_{n_1} + \hat{\beta}_{n_1})/2)(I(U_1 \le d^0 + s/n_2^{1+\gamma}) - I(U_1 \le d^0))}|\eta_{n_1,1} = \eta_1, \\
&\qquad\qquad\qquad\qquad\qquad\qquad\qquad\qquad \eta_{n_1,2} = \eta_2, \xi_{n_1} = \xi]^{n_2}.
\end{aligned}$$

Let $\epsilon > 0$ be pre-assigned. By Proposition 1, we can find $L > 0$ such that for all sufficiently large $n$, $Z_{n_1}([-L, L]^3) \ge 1 - \epsilon/3$. Using the fact that characteristic functions are bounded by 1, it follows immediately that for all

$n \geq N_0$ (for some $N_0$),

$$|Q_{n_2,s}(t) - \phi_s(t)| \leq \int_{[-L,L]^3} |Q^\star_{n_2,s}(t,\eta_1,\eta_2,\xi) - \phi_s(t)| \, dZ_{n_1}(\eta_1,\eta_2,\xi) + 2\epsilon/3$$

$$\leq \sup_{(\eta_1,\eta_2,\xi) \in [-L,L]^3} |Q^\star_{n_2,s}(t,\eta_1,\eta_2,\xi) - \phi_s(t)| + 2\epsilon/3.$$

For this fixed $L$, we now show that for all sufficiently large $n$

$$D_n \equiv \sup_{(\eta_1,\eta_2,\xi) \in [-L,L]^3} |Q^\star_{n_2,s}(t,\eta_1,\eta_2,\xi) - \phi_s(t)| \leq \epsilon/3,$$

whence it follows that eventually $|Q_{n_2,s}(t) - \phi_s(t)| \leq \epsilon$. To show the uniform convergence of $Q^\star_{n_2,s}(t,\eta_1,\eta_2,\xi)$ to $\phi_s(t)$ over the compact rectangle $[-L,L]^3$ we proceed as follows.

For given $L$ and $C$, it is the case that for all sufficiently large $n$, for any $\xi \in [-L,L]$ and any $0 < s < C$,

$$d^0 + \xi/n_1 - Kn_1^{-\gamma} < d^0 < d^0 + s/n_2^{1+\gamma} < d^0 + \xi/n_1 + Kn_1^{-\gamma}.$$

Let $P_{n_2,d}(s) \equiv \Pr(d^0 \leq U_1 \leq d^0 + s/n_2^{1+\gamma} | \hat{d}_{n_1} = d)$. Consider the conditional characteristic function $Q^\star_{n_2,s}(t,\eta_1,\eta_2,\xi)$, for $(\eta_1,\eta_2,\xi) \in [-L,L]^3$. It follows from the above display that for all sufficiently large $n$ (depending only on $L$ and $C$),

$$Q^\star_{n_2}(t,\eta_1,\eta_2,\xi)$$

$$= \Big[[1 - P_{n_2,d^0+\xi/n_1}(s)]$$

$$+ \int_{d^0}^{d^0+s/n_2^{1+\gamma}} E[e^{it(W_1 - (\hat{\alpha}_{n_1} + \hat{\beta}_{n_1})/2)} | U_1 = u] p_{U_1}(u) \, du\Big]^{n_2}$$

$$\left(\text{where } \hat{\alpha}_{n_1} = \alpha_0 + \frac{\eta_1}{\sqrt{n_1}}, \hat{\beta}_{n_1} = \beta_0 + \frac{\eta_2}{\sqrt{n_1}}\right)$$

$$= \Big[[1 - P_{n_2,d^0+\xi/n_1}(s)]$$

$$+ \frac{n_1^\gamma}{2K} \int_{d^0}^{d^0+s/n_2^{1+\gamma}} E[e^{it(W_1 - (\hat{\alpha}_{n_1} + \hat{\beta}_{n_1})/2)} | U_1 = u] \, du\Big]^{n_2}$$

$$= \Big[1 - \frac{1}{n_2} \frac{s}{2K} \left(\frac{\lambda}{1-\lambda}\right)^\gamma$$

$$+ \frac{n_1^\gamma}{2K} \int_{d^0}^{d^0+s/n_2^{1+\gamma}} E[e^{it(W_1 - (\hat{\alpha}_{n_1} + \hat{\beta}_{n_1})/2)} | U_1 = u] \, du\Big]^{n_2}$$

$$= \Big[1 - \frac{1}{n_2} \frac{s}{2K} \left(\frac{\lambda}{1-\lambda}\right)^\gamma + \frac{n_1^\gamma}{2K}$$

$$\times \int_0^s E\left[\exp\left\{it\left(W_1 - \frac{\alpha_0 + \eta_1/\sqrt{n_1} + \beta_0 + \eta_2/\sqrt{n_1}}{2}\right)\right\}\Big|\right.$$
$$\left.U_1 = d^0 + \frac{v}{n_2^{1+\gamma}}\right]\frac{1}{n_2^{1+\gamma}}\,dv\right]^{n_2}$$

$$= \left[1 - \frac{1}{n_2}\frac{s}{2K}\left(\frac{\lambda}{1-\lambda}\right)^\gamma\right.$$
$$+ \frac{1}{2Kn_2}\left(\frac{\lambda}{1-\lambda}\right)^\gamma e^{-it(\eta_1+\eta_2)/(2\sqrt{n_1})}$$
$$\left.\times \int_0^s E\left[e^{it(W_1-(\alpha_0+\beta_0)/2)}|U = d^0 + \frac{v}{n_2^{1+\gamma}}\right]dv\right]^{n_2}$$
$$= \left[1 - \frac{1}{n_2}\frac{s}{2K}\left(\frac{\lambda}{1-\lambda}\right)^\gamma(1 - B_{n_1,n_2,\eta_1,\eta_2}(s))\right]^{n_2},$$

where

$$B_{n_1,n_2,\eta_1,\eta_2}(s) = \frac{1}{s}e^{-it(\eta_1+\eta_2)/(2\sqrt{n_1})}\int_0^s E\left[e^{it(W_1-(\alpha_0+\beta_0)/2)}|U_1 = d^0 + \frac{v}{n_2^{1+\gamma}}\right]dv$$

and

$$D_n = \sup_{(\eta_1,\eta_2)\in[-L,L]^2}\left|\left[1 - \frac{1}{n_2}\frac{s}{2K}\left(\frac{\lambda}{1-\lambda}\right)^\gamma(1 - B_{n_1,n_2,\eta_1,\eta_2}(s))\right]^{n_2} - \phi_s(t)\right|.$$

Let $z_n(\eta_1, \eta_2) = -\frac{s}{2K}(\frac{\lambda}{1-\lambda})^\gamma(1 - B_{n_1,n_2,\eta_1,\eta_2}(s))$. It is easy to see that

$$\tilde{D}_n \equiv \sup_{(\eta_1,\eta_2)\in[-L,L]^2}|z_n(\eta_1,\eta_2) - z_0| \to 0,$$

where $z_0 = -(s/2K)(\lambda/(1-\lambda))^\gamma(1 - L^+(t))$. Consider now

$$D_n = \sup_{(\eta_1,\eta_2)\in[-L,L]^2}\left|\left(1 + \frac{1}{n_2}z_n(\eta_1,\eta_2)\right)^{n_2} - e^{z_0}\right|.$$

This is dominated by $I_n + II_n$ where

$$I_n = \left|\left(1 + \frac{1}{n_2}z_0\right)^{n_2} - e^{z_0}\right| \to 0$$

and

$$II_n = \sup_{(\eta_1,\eta_2)\in[-L,L]^2}\left|\left(1 + \frac{1}{n_2}z_n(\eta_1,\eta_2)\right)^{n_2} - \left(1 + \frac{1}{n_2}z_0\right)^{n_2}\right|.$$

Since $\tilde{D}_n$ goes to 0, for all sufficiently large $n$, $|z_0| \vee (\sup_{\eta_1,\eta_2\in[-L,L]^2}|z_n(\eta_1,\eta_2)|)$ is bounded by a constant, say $M$. Straightforward algebra shows that for all

sufficiently large $n$,

$$II_n \leq \left(\sup_{(\eta_1,\eta_2)\in[-L,L]^2} |z_n(\eta_1,\eta_2) - z_0|\right) \left(\sum_{j=1}^{n_2} \binom{n_2}{j} \frac{jM^{j-1}}{n_2^j}\right)$$

$$= \left(\sup_{(\eta_1,\eta_2)\in[-L,L]^2} |z_n(\eta_1,\eta_2) - z_0|\right) \left(1 + \frac{M}{n_2}\right)^{n_2-1} \to 0.$$

Thus, $D_n \to 0$ and the uniform convergence of $Q^{\star}_{n_2}(t,\eta_1,\eta_2,\xi)$ to $\phi_s(t) = e^{z_0}$ on $[-L,L]^3$ is established.

The next step is to establish the weak convergence of the finite-dimensional distributions of $(\mathbb{M}_{n_2}, \mathbb{J}_{n_2})$ to those of $(\mathbb{M}, \mathbb{J})$. For convenience, we restrict ourselves only to the set $[0, C]$. Let $J$ be a positive integer and consider $0 = s_0 < s_1 < s_2 < \cdots < s_J \leq C$. Let $c_1, c_2, \ldots, c_J$ and $d_1, d_2, \ldots, d_J$ be constants. We want to show that

$$A_n \equiv E(e^{it\sum_{j\leq J}(c_j(\mathbb{M}^+_{n_2}(s_j) - \mathbb{M}^+_{n_2}(s_{j-1})) + d_j(\mathbb{J}^+_{n_2}(s_j) - \mathbb{J}^+_{n_2}(s_{j-1})))})$$

$$\to A \equiv E(e^{it\sum_{j\leq J}(c_j(\mathbb{M}(s_j) - \mathbb{M}(s_{j-1})) + d_j(\mathbb{J}(s_j) - \mathbb{J}(s_{j-1})))})$$

for any vector of constants $(c_1, c_2, \ldots, c_j, d_1, d_2, \ldots, d_j)$. By the Cramer–Wold device it follows that

$$(\{\mathbb{M}_{n_2}(s_i) - \mathbb{M}_{n_2}(s_{i-1})\}_{i=1}^J, \{\mathbb{J}_{n_2}(s_i) - \mathbb{J}_{n_2}(s_{i-1})\}_{i=1}^J)$$

$$\xrightarrow{d} (\{\mathbb{M}(s_i) - \mathbb{M}(s_{i-1})\}_{i=1}^J, \{\mathbb{J}(s_i) - \mathbb{J}(s_{i-1})\}_{i=1}^J),$$

establishing the claim. As before, $A_n = \int K^{\star}_{n_2}(t,\eta_1,\eta_2,\xi)\, dZ_{n_1}(\eta_1,\eta_2,\xi)$, where

$$K^{\star}_{n_2}(t,\eta_1,\eta_2,\xi) = E[e^{it\sum_{j\leq J}(c_j(\mathbb{M}^+_{n_2}(s_j) - \mathbb{M}^+_{n_2}(s_{j-1})) + d_j(\mathbb{J}^+_{n_2}(s_j) - \mathbb{J}^+_{n_2}(s_{j-1})))}|$$

$$\eta_{n_1,1} = \eta_1, \eta_{n_1,2} = \eta_2, \xi_{n_1} = \xi].$$

Proceeding as before, the convergence of $A_n$ to $A$ follows if we establish the uniform convergence of $K^{\star}_{n_2}(t,\eta_1,\eta_2,\xi)$ to $A$ on a compact rectangle of the form $[-L,L]^3$. This is achieved by using arguments similar to those involved in the proof of the convergence of the one-dimensional marginals above. The details of the algebra are available in the Appendix of Lan, Banerjee and Michailidis (2007). The derivation there can be extended readily to allow for $s_j$'s that can also be negative but that has been avoided as that extension involves no new ideas but becomes more cumbersome.

We finally show that the process $\mathbb{M}_{n_2}$ restricted to $[-C, C]$ is tight. We know that $(\hat{\alpha}_{n_1}, \hat{\beta}_{n_1}, \hat{d}_{n_1}) \longrightarrow_p (\alpha_0, \beta_0, d^0)$ and $n_1(\hat{d}_{n_1} - d^0) = O_p(1)$. Let

$$\Omega_n = \left\{|\hat{\alpha}_{n_1} - \alpha_0| \leq \Delta, |\hat{\beta}_{n_1} - \beta_0| \leq \Delta,\right.$$

$$\left.\hat{d}_{n_1} - \frac{K}{n_1^\gamma} < d^0 - \frac{C}{n_2^{1+\gamma}} < d^0 + \frac{C}{n_2^{1+\gamma}} < \hat{d}_{n_1} + \frac{K}{n_1^\gamma}\right\}.$$

Clearly, $P(\Omega_n) \longrightarrow 1$. The event $\Omega_n$ can be written as $(\hat{\alpha}_{n_1}, \hat{\beta}_{n_1}, \hat{d}_{n_1}) \in R_n$, where $H_n(R_n) \longrightarrow 1$, $H_n$ being the joint distribution of $(\hat{\alpha}_{n_1}, \hat{\beta}_{n_1}, \hat{d}_{n_1})$. Note that $\mathbb{M}_{n_2}\mathbf{1}(\Omega_n)$ is also a process in $D(\mathbb{R})$. We verify tightness of $\mathbb{M}_{n_2}\mathbf{1}(\Omega_n)$ restricted to $[-C, C]$. To this end, we verify (the analogue of) condition (13.14) on page 143 of Billingsley (1999), with $\beta = 1/2$ and $\alpha = 1$. Once again, let $0 \leq s_1 \leq s \leq s_2 \leq C$,

$$E[|\mathbb{M}_{n_2}^+(s) - \mathbb{M}_{n_2}^+(s_1)| \cdot |\mathbb{M}_{n_2}^+(s_2) - \mathbb{M}_{n_2}^+(s)|\mathbf{1}(\Omega_n)]$$

$$= \int_{R_n} E[|\mathbb{M}_{n_2}^+(s) - \mathbb{M}_{n_2}^+(s_1)| \cdot |\mathbb{M}_{n_2}^+(s_2) - \mathbb{M}_{n_2}^+(s)|$$

$$\times |(\hat{\alpha}_{n_1}, \hat{\beta}_{n_1}, \hat{d}_{n_1}) = (\alpha, \beta, d)] \, dH_n(\alpha, \beta, d),$$

where

$$E[|\mathbb{M}_{n_2}^+(s) - \mathbb{M}_{n_2}^+(s_1)| \cdot |\mathbb{M}_{n_2}^+(s_2) - \mathbb{M}_{n_2}^+(s)||(\hat{\alpha}_{n_1}, \hat{\beta}_{n_1}, \hat{d}_{n_1}) = (\alpha, \beta, d)]$$

$$= E_{\alpha,\beta,d}\left[\left|\sum_{i=1}^{n_2}\left(W_i - \frac{\alpha+\beta}{2}\right)\left(I\left(U_i \leq d^0 + \frac{s}{n_2^{1+\gamma}}\right)\right.\right.\right.$$

$$\left.\left.- I\left(U_i \leq d^0 + \frac{s_1}{n_2^{1+\gamma}}\right)\right)\right|$$

$$\times \left|\sum_{i=1}^{n_2}\left(W_i - \frac{\alpha+\beta}{2}\right)\left(I\left(U_i \leq d^0 + \frac{s_2}{n_2^{1+\gamma}}\right)\right.\right.$$

$$\left.\left.\left.- I\left(U_i \leq d^0 + \frac{s}{n_2^{1+\gamma}}\right)\right)\right|\right]$$

$$= E_{\alpha,\beta,d}\left[\left|\sum_{i \neq j} I\left(d^0 + \frac{s_1}{n_2^{1+\gamma}} \leq U_i < d^0 + \frac{s}{n_2^{1+\gamma}}\right)\right.\right.$$

$$\times I\left(d^0 + \frac{s}{n_2^{1+\gamma}} \leq U_j < d^0 + \frac{s_2}{n_2^{1+\gamma}}\right)$$

$$\left.\left.\times \left(W_i - \frac{\alpha+\beta}{2}\right)\left(W_j - \frac{\alpha+\beta}{2}\right)\right|\right]$$

$$\leq \sum_{i \neq j} E_{\alpha,\beta,d}\left[\left|I\left(d^0 + \frac{s_1}{n_2^{1+\gamma}} \leq U_i < d^0 + \frac{s}{n_2^{1+\gamma}}\right)\right.\right.$$

$$\times I\left(d^0 + \frac{s}{n_2^{1+\gamma}} \leq U_j < d^0 + \frac{s_2}{n_2^{1+\gamma}}\right)$$

$$\left.\left.\times \left(W_i - \frac{\alpha+\beta}{2}\right)\left(W_j - \frac{\alpha+\beta}{2}\right)\right|\right]$$

$$= \sum_{i \neq j} \sup_{t \in [d^0 + s_1 n_2^{-(1+\gamma)}, d^0 + s n_2^{-(1+\gamma)}]} E_{\alpha,\beta,d}\left[\left|W_i - \frac{\alpha+\beta}{2}\right|\Big|U_i = t\right] \frac{n_1^\gamma}{2K} \frac{s - s_1}{n_2^{1+\gamma}}$$

$$\times \sup_{t \in [d^0 + s n_2^{-(1+\gamma)}, d^0 + s_2 n_2^{-(1+\gamma)}]} E_{\alpha,\beta,d}\left[\left|W_j - \frac{\alpha+\beta}{2}\right|\Big|U_j = t\right] \frac{n_1^\gamma}{2K} \frac{s_2 - s}{n_2^{1+\gamma}}$$

$$\leq \sum_{i \neq j} H_1\left(\frac{s - s_1}{2K n_2}\left(\frac{\lambda}{1-\lambda}\right)^\gamma\right) H_2\left(\frac{s_2 - s}{2K n_2}\left(\frac{\lambda}{1-\lambda}\right)^\gamma\right)$$

for any $(\alpha, \beta, d) \in R_n$

$$= H_1 H_2 \sum_{i \neq j} \frac{(s - s_1)(s_2 - s)}{n_2^2 4k^2}\left(\frac{\lambda}{1-\lambda}\right)^{2\gamma} \leq c^*(s_2 - s_1)^2.$$

It follows that

$$E[|\mathbb{M}_{n_2}^+(s) - \mathbb{M}_{n_2}^+(s_1)| \cdot |\mathbb{M}_{n_2}^+(s_2) - \mathbb{M}_{n_2}^+(s)| \mathbf{1}(\Omega_n)]$$
$$\leq H_n(R_n) c^*(s_2 - s_1)^2 \leq c^*(s_2 - s_1)^2,$$

which establishes tightness of $\mathbb{M}_{n_2} \mathbf{1}(\Omega_n)$.

Then, given any $\epsilon > 0$, $\forall n > N_1$,

$$\text{Prob}[\omega : \mathbb{M}_{n_2} \mathbf{1}(\Omega_n)(\omega) \in \mathcal{K}] \geq 1 - \varepsilon,$$

where $\mathcal{K}$ is a compact set. But $\text{Prob}[\omega : \omega \in \Omega_n] \geq 1 - \varepsilon$ eventually. Therefore, eventually

$$\text{Prob}[\omega \in \Omega_n \text{ and } \mathbb{M}_{n_2} \mathbf{1}(\Omega_n) \in \mathcal{K}] \geq 1 - 2\varepsilon$$

and consequently $\text{Prob}[\mathbb{M}_{n_2} \in \mathcal{K}] \geq 1 - 2\varepsilon$. This establishes the tightness of $\mathbb{M}_{n_2}$ in the space of right-continuous left limits endowed functions on $[-C, C]$. Similarly, the tightness of $\mathbb{J}_{n_2}$ can be established. This completes the verification of marginal tightness and therefore joint tightness.

Before embarking on the proof of Theorem 1, we need some auxiliary lemmas. We first state these below.

LEMMA A.1. *Suppose that $X_1, X_2, \ldots, X_n$ are i.i.d. random elements assuming values in a space $\mathcal{X}$. Let $\mathcal{F}$ be a class of functions with domain $\mathcal{X}$ and range in $[0,1]$ with finite VC dimension $V(\mathcal{F})$ and set $V = 2(V(\mathcal{F}) - 1)$. Denoting by $\mathbb{P}_n$ the empirical measure corresponding to the sample and by $P$ the distribution of $X_1$, we have*

$$\text{Pr}^\star(\|\sqrt{n}(\mathbb{P}_n - P)\|_{\mathcal{F}} \geq \lambda) \leq \left(\frac{D\lambda}{\sqrt{V}}\right)^V \exp(-2\lambda^2).$$

This lemma is adapted from Talagrand (1994).

LEMMA A.2. *Let $\tilde{U}_1, \tilde{U}_2, \ldots, \tilde{U}_n$ be i.i.d. random variables following the uniform distribution on $(0,1)$. Let $\tilde{\epsilon}_1, \tilde{\epsilon}_2, \ldots, \tilde{\epsilon}_n$ be i.i.d. mean 0 random variables with finite variance $\sigma^2$ that are independent of the $\tilde{U}_i$'s. Let $\beta_n(s) = \sum_{i=1}^n \tilde{\epsilon}_i 1(\tilde{U}_i \leq s)$. Then for any $0 < \alpha < \beta < 1$, we have*

$$\Pr\left(\sup_{\alpha \leq s \leq \beta} \frac{|\beta_n(s)|}{s} \geq \lambda\right) \leq (\alpha^{-1} - \beta^{-1})\lambda^{-2}\sigma^2.$$

The proof of this lemma follows the proof of Theorem 1. □

PROOF OF THEOREM 1. For simplicity and ease of exposition, in what follows we assume that $n$ points are used at the first-stage to compute estimates $\hat{\alpha}_n, \hat{\beta}_n, \hat{d}_{n_1}$ of the three parameters of interest. At the second stage $n$ i.i.d. $U_i$'s are sampled from the uniform distribution on $\tilde{D}_n \equiv [\hat{d}_{n_1} - Kn^{-\gamma}, \hat{d}_{n_1} + Kn^{-\gamma}]$ and the updated estimate of $d^0$ is computed as

$$\hat{d}_{n_2} = \arg\min_{u \in \tilde{D}_n} \frac{1}{n} \sum_{i=1}^n [W_i - \hat{\alpha}_n 1(U_i \leq u) - \hat{\beta}_n 1(U_i > u)]^2 \equiv \arg\min_{u \in \tilde{D}_n} SS(u).$$

In the above display $W_i = f(U_i) + \epsilon_i$ where $\epsilon_i$'s are i.i.d. error variables. Working under this more restrictive setting (of equal allocation of points at each stage) *does not compromise* the complexity of the arguments involved. Finally, recall that by our assumption, $E[e^{C|\epsilon_1|}]$ is finite, for some $C > 0$.

Before proceeding further, a word about the definition of arg min in the above display. The function $SS$ is a right-continuous function endowed with left limits. For this derivation, we take the arg min to be the smallest $u$ in $\tilde{D}_n$ for which $\min(SS(u), SS(u-)) = \inf_{x \in \tilde{D}_n} SS(x)$.

Denote by $G_n$ the distribution of $(\hat{\alpha}_n, \hat{\beta}_n, \hat{d}_{n_1})$. Now, given $\epsilon > 0$, find $L$ so large that for all sufficiently large $n$, say $n \geq N_0$,

$$(\hat{\alpha}_n, \hat{\beta}_n, \hat{d}_{n1}) \in [\alpha_0 - L/\sqrt{n}, \alpha_0 + L/\sqrt{n}] \times [\beta_0 - L/\sqrt{n}, \beta_0 + L/\sqrt{n}]$$
$$\times [d^0 - L/n, d^0 + L/n]$$

with probability greater than $1 - \epsilon$. Denote the region on the right-hand side of the above display by $R_n$. Then, for all $n \geq N_0$,

$$\Pr(n^{1+\gamma}|\hat{d}_{n_2} - d^0| > a)$$
$$\leq \int_{R_n} \Pr(n^{1+\gamma}|\hat{d}_{n_2} - d^0| > a | \hat{\alpha}_n = \alpha, \hat{\beta}_n = \beta, \hat{d}_{n_1} = t) \, dG_n(\alpha, \beta, t) + \epsilon,$$

which is dominated by $\sup_{(\alpha,\beta,t) \in R_n} \Pr_{t,\alpha,\beta}(n^{1+\gamma}|\hat{d}_{n_2} - d^0| > a) + \epsilon$. By making $a$ large, we will show that for all sufficiently large $n$ (say $n > N_1 > N_0$), the supremum is bounded by $\epsilon$. This will complete the proof.

First note that since $N_0$ is chosen to be sufficiently large, whenever $n \geq N_0$ and $t \in [d^0 - L/n, d^0 + L/n]$, it is the case that $t - Kn^{-\gamma} < d^0 - Kn^{-\gamma}/2 <$ $d^0 + Kn^{-\gamma}/2 < t + Kn^{-\gamma}]$. Let $\tilde{P}_n$ be the distribution of the pair $(W_1, U_1)$ generated at stage two under first-stage parameters $(\alpha, \beta, t)$ and $\tilde{\mathbb{P}}_n$ is the empirical measure corresponding to $n$ i.i.d. observations from $\tilde{P}_n$. It is not difficult to see that $\hat{d}_{n_2} \equiv \arg\min_{d \in [t-Kn^{-\gamma}, t+Kn^{-\gamma}]} \tilde{\mathbb{M}}_n(d)$ where $\tilde{\mathbb{M}}_n(d) = \tilde{\mathbb{P}}_n[(w - (\alpha + \beta)/2)(1(u \leq d) - 1(u \leq d^0))]$, while, $d = \arg\min_{d \in [t-Kn^{-\gamma}, t+Kn^{-\gamma}]} \tilde{M}_n(d)$ where $\tilde{M}_n(d) = \tilde{P}_n[(w - (\alpha + \beta)/2)(1(u \leq d) - 1(u \leq d^0))]$. We have

$$\tilde{M}_n(d) = \{|\beta_0 - (\alpha + \beta)/2||d - d^0|1(d \geq d^0)$$
$$+ |\alpha_0 - (\alpha + \beta)/2||d - d^0|1(d < d^0)\}(n^\gamma/2K).$$

Now, for $0 < r \leq K/2$, set $a(r) = \min\{\tilde{M}_n(d) : |d - d^0| \geq rn^{-\gamma}\}$. Then $a(r) = \min(|\beta_0 - (\alpha + \beta)/2)|, |\alpha_0 - (\alpha + \beta)/2)|)r/2K$ and let $b(r) = (a(r) - \tilde{M}_n(d^0))/3 = a(r)/3$. Now, for all $n \geq N_0$, for $\alpha, \beta$ in the region under consideration, $b(r)$ is readily seen to be uniformly bounded below by $\kappa r$ for some constant $\kappa$ depending only on $\alpha_0, \beta_0, K, N_0$. We then have:

$$(13) \quad \sup_{d \in [t-Kn^{-\gamma}, t+Kn^{-\gamma}]} |\tilde{\mathbb{M}}_n(d) - \tilde{M}_n(d)| \leq b(r) \quad \Rightarrow \quad |\hat{d}_{n_2} - d^0| \leq rn^{-\gamma}.$$

To prove this, assume that the inequality on the left-hand side of the above display holds and consider $d \in [t - Kn^{-\gamma}, t + Kn^{-\gamma}]$ with $|d - d^0| > rn^{-\gamma}$. Then, $\tilde{\mathbb{M}}_n(d) \geq \tilde{M}_n(d) - b(r) \geq a(r) - b(r)$ and $\tilde{\mathbb{M}}_n(d^o) \leq \tilde{M}_n(d^0) + b(r)$ jointly imply that $\tilde{\mathbb{M}}_n(d) - \tilde{\mathbb{M}}_n(d^0) \geq a(r) - b(r) - \tilde{M}_n(d^0) - b(r) = b(r) > 0$. Hence, $\tilde{\mathbb{M}}_n(d) > \tilde{\mathbb{M}}_n(d^0)$.

Now, since $\hat{d}_{n_2}$ is the smallest $d \in \tilde{D}_n$ for which $\tilde{\mathbb{M}}_n(d) \wedge \tilde{\mathbb{M}}_n(d-) = \inf_{x \in \tilde{D}_n} \tilde{\mathbb{M}}_n(x)$ and $\tilde{\mathbb{M}}_n$ is a (right-continuous left limits endowed) piecewise constant function with finitely many flat stretches, it is easy to see that $\tilde{\mathbb{M}}_n(\hat{d}_{n_2}) = \inf_{x \in \tilde{D}_n} \tilde{\mathbb{M}}_n(x)$. Therefore, $\tilde{\mathbb{M}}_n(\hat{d}_{n_2}) \leq \tilde{\mathbb{M}}_n(d^0)$, showing that $|\hat{d}_{n_2} - d^0| \leq rn^{-\gamma}$ in view of the last display above.

Now, consider $\Pr_{\alpha, \beta, t}(|\hat{d}_{n_2} - d^0| > rn^{-\gamma})$, which is

$$\Pr_{\alpha, \beta, t}(|\hat{d}_{n_2} - d^0| > rn^{-\gamma})$$
$$(14) \quad \leq \Pr_{\alpha, \beta, t}(rn^{-\gamma} < |\hat{d}_{n_2} - d^0| \leq \delta n^{-\gamma})$$
$$+ \Pr_{\alpha, \beta, t}(|\hat{d}_{n_2} - d^0| > \delta n^{-\gamma})$$
$$(15) \quad \equiv P_n(\alpha, \beta, t) + Q_n(\alpha, \beta, t),$$

where $\delta$ (is sufficiently small, say less than $K/3$) does not depend on $t, \alpha, \beta$. We deal with $Q_n(\alpha, \beta, t)$ later. We first consider $P_n(\alpha, \beta, t) \equiv \Pr_{\alpha, \beta, t}(rn^{-\gamma} < |\hat{d}_{n_2} - d^0| \leq \delta n^{-\gamma})$. Since,

$$\{rn^{-\gamma} < |\hat{d}_{n_2} - d^0| \leq \delta n^{-\gamma}\} \subseteq \left[ \bigcup_{d^0 + rn^{-\gamma} < d \leq d^0 + \delta n^{-\gamma}} \{\tilde{\mathbb{M}}_n(d) \leq \tilde{\mathbb{M}}_n(d^0)\} \right]$$

$$\cup \left[\bigcup_{d^0-\delta n^{-\gamma}\leq d<d^0-rn^{-\gamma}} \{\tilde{\mathbb{M}}_n(d)\leq \tilde{\mathbb{M}}_n(d^0)\}\right],$$

we conclude that

$$P_n(\alpha,\beta,t) \leq P_{n,1}(\alpha,\beta,t) + P_{n,2}(\alpha,\beta,t)$$

$$\equiv \Pr_{\alpha,\beta,t}\left(\bigcup_{d^0+rn^{-\gamma}<d\leq d^0+\delta n^{-\gamma}} \{\tilde{\mathbb{M}}_n(d^0) - \tilde{\mathbb{M}}_n(d) \geq 0\}\right)$$

$$+ \Pr_{\alpha,\beta,t}\left(\bigcup_{d^0-\delta n^{-\gamma}\leq d<d^0-rn^{-\gamma}} \{\tilde{\mathbb{M}}_n(d^0) - \tilde{\mathbb{M}}_n(d) \geq 0\}\right).$$

We first construct an upper bound on $\sup_{(\alpha,\beta,t)\in R_n} P_{n,1}(\alpha,\beta,t)$. For any $d \in (d^0 + rn^{-\gamma}, d^0 + \delta n^{-\gamma}]$ we have

$$\tilde{\mathbb{M}}_n(d^0) - \tilde{\mathbb{M}}_n(d)$$

$$= (\tilde{\mathbb{M}}_n(d^0) - \tilde{M}_n(d^0)) - (\tilde{\mathbb{M}}_n(d) - \tilde{M}_n(d)) - (\tilde{M}_n(d) - \tilde{M}_n(d^0))$$

$$= -(\tilde{\mathbb{M}}_n(d) - \tilde{M}_n(d)) - \left|\beta_0 - \frac{\alpha+\beta}{2}\right|\frac{n^\gamma}{2K}|d-d^0|.$$

Hence $0 \leq \tilde{\mathbb{M}}_n(d^0) - \tilde{\mathbb{M}}_n(d) \Rightarrow (2K)^{-1}|\beta_0 - \frac{\alpha+\beta}{2}| \leq \frac{-(\tilde{\mathbb{M}}_n(d)-\tilde{M}_n(d))}{n^\gamma|d-d^0|}$. Now, for all $(\alpha,\beta,t) \in R_n$ (with $n \geq N_0$), $|\beta_0 - \frac{\alpha+\beta}{2}|$ is bounded below by some constant $B$, whence it follows that $0 \leq \tilde{\mathbb{M}}_n(d^0) - \tilde{\mathbb{M}}_n(d) \Rightarrow \frac{|\tilde{\mathbb{M}}_n(d)-\tilde{M}_n(d)|}{n^\gamma|d-d^0|} \geq \frac{B}{2K}$. Thus,

$$\bigcup_{d^0-\delta n^{-\gamma}<d\leq d^0-rn^{-\gamma}} \{\tilde{\mathbb{M}}_n(d^0) - \tilde{\mathbb{M}}_n(d) \geq 0\}$$

$$\subset \left\{\sup_{d^0+rn^{-\gamma}<d\leq d^0+\delta n^{-\gamma}} \frac{|\tilde{\mathbb{M}}_n(d)-\tilde{M}_n(d)|}{n^\gamma|d-d^0|} \geq \tilde{B}\right\},$$

where $\tilde{B} = B/2K$. We thus have

$$P_{n,1}(\alpha,\beta,t) \leq \Pr_{\alpha,\beta,t}\left[\sup_{d^0+rn^{-\gamma}<d\leq d^0+\delta n^{-\gamma}} \frac{|\tilde{\mathbb{M}}_n(d)-\tilde{M}_n(d)|}{n^\gamma|d-d^0|} \geq \tilde{B}\right]$$

$$= \Pr_{\alpha,\beta,t}\left(\sup_{d^0+rn^{-\gamma}<d\leq d^0+\delta n^{-\gamma}} \frac{|\sqrt{n}(\tilde{\mathbb{P}}_n - \tilde{P}_n)f_{d,\alpha,\beta}(u,w)|}{|d-d^0|n^\gamma} \geq \sqrt{n}\tilde{B}\right),$$

where

$$f_{d,\alpha,\beta}(u,w) = (w - (\alpha+\beta)/2)(1(u \leq d) - 1(u \leq d^0)).$$

Using the fact that for $d > d^0$, $(W_j - (\alpha+\beta)/2)(1(U_j \leq d) - 1(U_j \leq d^0)) = (\beta_0 - (\alpha+\beta)/2)(1(U_j \leq d) - 1(U_j \leq d^0)) + \epsilon_j(1(U_j \leq d) - 1(U_j \leq d^0))$, this upper bound on $P_{n,1}(\alpha,\beta,t)$ is easily seen to be dominated by $I_n + II_n$ where

$$I_n = \Pr\nolimits_{\alpha,\beta,t}\bigg(|\beta_0 - (\alpha+\beta)/2|$$

$$\times \sup_{r<s\leq\delta} \frac{|\sqrt{n}(\tilde{\mathbb{P}}_n - \tilde{P}_n)(1(u \leq d^0 + sn^{-\gamma}) - 1(u \leq d^0))|}{s}$$

$$\geq \sqrt{n}\tilde{B}/2\bigg),$$

which in turn is dominated by

$$I'_n = \Pr\nolimits_{\alpha,\beta,t}\bigg(\sup_{r<s\leq\delta} \frac{|\sqrt{n}(\tilde{\mathbb{P}}_n - \tilde{P}_n)(1(u \leq d^0 + sn^{-\gamma}) - 1(u \leq d^0))|}{s} \geq \frac{\sqrt{n}\tilde{B}}{2B'}\bigg),$$

where, for $n \geq N_0$ and $(\alpha,\beta,t) \in R_n$, $B'$ is an upper bound on $|\beta_0 - (\alpha+\beta)/2|$, while

$$II_n = \Pr\nolimits_t\bigg(\sup_{r<s\leq\delta} \frac{|(1/\sqrt{n})\sum_{i=1}^n \epsilon_i(1(U_i \leq d^0 + sn^{-\gamma}) - 1(U_i \leq d^0))|}{s}$$

$$\geq \sqrt{n}\tilde{B}/2\bigg).$$

Since the $U_i$'s are i.i.d. Uniform on $[t - Kn^{-\gamma}, t + Kn^{-\gamma}]$, it is easy to see that $I'_n$ is simply

$$\Pr\bigg(\sup_{r<s\leq\delta} \frac{|\sqrt{n}(\mathbb{Q}_n - Q)(1(\tilde{w} \leq s))|}{s} \geq \frac{\sqrt{n}\tilde{B}}{2B'}\bigg),$$

where $\tilde{W}_1, \tilde{W}_2, \ldots, \tilde{W}_n$ are i.i.d. Unif$[0, 2K]$, $\mathbb{Q}_n$ is the empirical measure of the $\tilde{W}_i$s and $Q$ is the distribution of $\tilde{W}_1$. In terms of $\tilde{U}_1, \tilde{U}_2, \ldots, \tilde{U}_n$, which are i.i.d. Unif$[0,1]$, this expression is simply:

$$\Pr\bigg(\sup_{r/2K<s\leq\delta/2K} \frac{|\sqrt{n}(\mathbb{P}_n - P)(1(\tilde{u} \leq s))|}{s} \geq \frac{\sqrt{n}2K\tilde{B}}{2B'}\bigg).$$

By Lemma A.3 of Ferger (2005), this is bounded above by a constant (that depends only on $\alpha_0, \beta_0, K, N_0$) times $1/rn$. Now, in terms of $\tilde{U}_1, \ldots, \tilde{U}_n$ and $\tilde{\epsilon}_1, \ldots, \tilde{\epsilon}_n$ (where the $\tilde{\epsilon}_i$'s are defined on the same probability space as the $\tilde{U}_i$'s, but independently of them, and are distributed like the $\epsilon_i$'s) $II_n$ is simply

$$\Pr\bigg(\sup_{r/2K<s\leq\delta/2K} \frac{|(1/\sqrt{n})\sum_{i=1}^n \tilde{\epsilon}_i(1(\tilde{U}_i \leq s))|}{s} \geq \frac{\sqrt{n}2K\tilde{B}}{2}\bigg)$$

and this, by Lemma A.2, is dominated up to a constant (that only depends on $\alpha_0, \beta_0, \sigma, K, N_0$) by $(1/rn)$. It follows that for some constant $C_0$, for all $n \geq N_0$, $\sup_{(\alpha,\beta,t)\in R_n} P_{n1}(\alpha,\beta,t) \leq C_0/rn$. A similar (uniform) bound works $P_{n2}(\alpha,\beta,t)$. It follows that $\sup_{(\alpha,\beta,t)\in R_n} P_n(\alpha,\beta,t) \leq C_0/rn$, at the expense of a larger constant $C_0$. Thus, from (15), we have

$$\sup_{(\alpha,\beta,t)\in R_n} \Pr_{\alpha,\beta,t}(|\hat{d}_{n_2} - d^0| > rn^{-\gamma}) \leq C_0(rn)^{-1} + \sup_{(\alpha,\beta,t)\in R_n} Q_n(\alpha,\beta,t).$$

To find a uniform upper bound on $Q_n(\alpha,\beta,t)$ note that, from (13), we have, for all $n > N_0$

$$\Pr_{\alpha,\beta,t}(|\hat{d}_{n_2} - d^0| > \delta n^{-\gamma})$$

$$\leq \Pr_{\alpha,\beta,t}\left(\sup_{d\in[t-Kn^{-\gamma},t+Kn^{-\gamma}]} |\tilde{\mathbb{M}}_n(d) - \tilde{M}_n(d)| > b(\delta)\right)$$

$$\leq \Pr_{\alpha,\beta,t}\left(\sup_{d\in[t-Kn^{-\gamma},t+Kn^{-\gamma}]} |\tilde{\mathbb{M}}_n(d) - \tilde{M}_n(d)| > \kappa\delta\right)$$

and it suffices to find a uniform upper bound for this last expression. But this is bounded by

$$\Pr_{\alpha,\beta,t}\left[\sup_{d\in[t-Kn^{-\gamma},t+Kn^{-\gamma}]} |\sqrt{n}(\tilde{\mathbb{P}}_n - \tilde{P}_n)(\mu(u) - (\alpha+\beta)/2)\right.$$

$$\times (1(u \leq d) - 1(u \leq d^0))| > \sqrt{n}\kappa\delta/2\bigg]$$

$$+ \Pr_{\alpha,\beta,t}\left[\sup_{d\in[t-Kn^{-\gamma},t+Kn^{-\gamma}]} \left|n^{-1}\sum_{i=1}^n \epsilon_i(1(U_i \leq d) - 1(U_i \leq d^0))\right| > \kappa\delta/2\right].$$

To tackle the first term, we invoke Lemma A.1. For $(\alpha,\beta,t) \in R_n$, the class $[(\mu(u) - (\alpha+\beta)/2)(1(u \leq d) - 1(u \leq d^0)) : d \in [t - Kn^{-\gamma}, t + Kn^{-\gamma}]]$ is a bounded VC class of functions (with the bound not depending on $\alpha,\beta,t$) and with finite VC dimension, say $\mathcal{V}$ (which does not depend on $\alpha,\beta,t$). Hence, we can apply Lemma A.1 to conclude that

$$\Pr_{\alpha,\beta,t}\left[\sup_{d\in[t-Kn^{-\gamma},t+Kn^{-\gamma}]} |\sqrt{n}(\tilde{\mathbb{P}}_n - \tilde{P}_n)(\mu(u) - (\alpha+\beta)/2)\right.$$

$$\times (1(u \leq d) - 1(u \leq d^0))| > \sqrt{n}\kappa\delta/2\bigg]$$

$$\leq \tilde{C}_1 \times (\sqrt{n}\kappa\delta)^{2(\mathcal{V}-1)} \exp(-\tilde{C}_2 n\kappa^2\delta^2),$$

where the constants $\tilde{C}_1$ and $\tilde{C}_2$ depend solely on the VC dimension and the upper bound on the functions. For all sufficiently large $n$, the right-hand side of the above display is less than $\epsilon/3$. To deal with the second term, we use the results on pages 132, 133 of Van de Geer (2000). We write the second term as:

$$\int \Pr_{\alpha,\beta,t}\left( \sup_{d \in [t-Kn^{-\gamma}, t+Kn^{-\gamma}]} \left| n^{-1} \sum_{i=1}^{n} \epsilon_i (1(u_i \leq d) - 1(u_i \leq d^0)) \right| > \kappa\delta/2 \right) dH_n(u_1, u_2, \ldots, u_n),$$

where $H_n$ is the joint distribution of $(U_1, U_2, \ldots, U_n)$. For each fixed $(u_1, u_2, \ldots, u_n)$, Corollary 8.8 of Van de Geer (2000) can be used to show that for $\delta$ sufficiently small and $n$ sufficiently large (where the thresholds do not depend on the $u_i$'s or $\alpha, \beta, t$),

$$\Pr_{\alpha,\beta,t}\left( \sup_{d \in [t-Kn^{-\gamma}, t+Kn^{-\gamma}]} \left| n^{-1} \sum_{i=1}^{n} \epsilon_i (1(u_i \leq d) - 1(u_i \leq d^0)) \right| > \kappa\delta/2 \right)$$
$$\leq \tilde{C} \exp(-\tilde{C}' n \delta^2)$$

for some constants $\tilde{C}$ and $\tilde{C}'$ that do not depend on $\alpha, \beta, t$ or the points $(u_1, u_2, \ldots, u_n)$. This implies that the second term can be made less than $\epsilon/3$ by choosing $n$ sufficiently large. It follows, that for all sufficiently large $n$ (say $n > N_1 > N_0$) and an appropriate choice of $\delta$, we have

$$\sup_{(\alpha,\beta,t) \in R_n} \Pr_{\alpha,\beta,t}(|\hat{d}_{n_2} - d^0| > r n^{-\gamma}) \leq C_0 (rn)^{-1} + 2\epsilon/3;$$

the first term on the right-hand side can be made less than $\epsilon/3$ by choosing $r = A/n$ where $A$ is large enough, showing that for all sufficiently large $n$, we can find $A$ large enough so that:

$$\sup_{(\alpha,\beta,t) \in R_n} \Pr_{\alpha,\beta,t}(n^{1+\gamma} |\hat{d}_{n_2} - d^0| > A) < \epsilon.$$

For the details as to how Corollary 8.8 of Van de Geer (2000) is applied in our setting, see the longer version of this proof in Lan, Banerjee and Michailidis (2007).  □

PROOF OF LEMMA A.2. Let $\beta_n(s) = \frac{1}{\sqrt{n}} \sum_{i=1}^{n} \tilde{\epsilon}_i 1(x_i \leq s)$. Let $\{s_k = \alpha + (\beta - \alpha) 2^{-m} : 0 \leq k \leq 2^m\}$, $m \in \mathbb{N}$ be a dyadic partition of $[\alpha, \beta]$. Consider

$P(\alpha, \beta, \lambda)$

$$= P\left( \sup_{\alpha \leq s \leq \beta} \frac{|\beta_n(s)|}{s} \geq \lambda \right)$$

$$= \lim_{m \to \infty} P\left(\max_{1 \leq k \leq 2^m} \frac{|\beta_n(s_k)|}{s_k} \geq \lambda\right)$$

$$= \lim_{m \to \infty} \int_{(x_1,\ldots,x_n) \in (0,1)^n} P\left(\max_{1 \leq k < 2^m} \left(\left|\frac{1}{\sqrt{n}} \sum_{i=1}^n \tilde{\epsilon}_i 1(x_i \leq s_k)\right|\right)\right.$$

$$\left.\times (s_k)^{-1} \geq \lambda\right) dx_1\, dx_2 \cdots dx_n.$$

For fixed $(x_1, x_2, \ldots, x_n) \in (0,1)^n$, set $M_k = \frac{1}{\sqrt{n}} \sum_{i=1}^n \tilde{\epsilon}_i 1(x_i \leq s_k), 0 \leq k \leq 2^m$. Define $X_k = M_k - M_{k-1}$ for $k \geq 1$. Then the $X_k$'s are independent random variables, each with mean 0 and finite variance and $M_k = X_1 + X_2 + \cdots + X_k$. Since $1/s_k$ is a decreasing sequence of constants, we can apply the Hajek–Renyi inequality [see, e.g., the Appendix of Lan, Banerjee and Michailidis (2007)] to conclude that

$$P\left(\max_{1 \leq k \leq 2^m} \frac{|M_k|}{s_k} \geq \lambda\right) \leq \frac{1}{\lambda^2} \sum_{1 \leq k \leq 2^m} s_k^{-2} E(M_k - M_{k-1})^2$$

$$= \frac{1}{\lambda^2} \sum_{1 \leq k \leq 2^m} s_k^{-2}(EM_k^2 - EM_{k-1}^2).$$

Now,

$$EM_k^2 = E\left[\left(\frac{1}{\sqrt{n}} \sum_{i=1}^n \tilde{\epsilon}_i 1(x_i \leq s_k)\right)^2\right] = \mathrm{Var}\left[\frac{1}{\sqrt{n}} \sum_{i=1}^n \tilde{\epsilon}_i 1(x_i \leq s_k)\right]^2$$

$$= \frac{1}{n}\sigma^2 \sum_{i=1}^n 1(x_i \leq s_k).$$

It follows that $EM_k^2 - EM_{k-1}^2 = \frac{\sigma^2}{n} \sum_{i=1}^n 1(s_{k-1} < x_i \leq s_k)$. Therefore

$$P(\alpha, \beta, \lambda) \leq \lim_{m \to \infty} \int_{(x_1,\ldots,x_n) \in (0,1)^n} \left\{\frac{1}{\lambda^2} \sum_{1 \leq k < 2^m} s_k^{-2} \frac{\sigma^2}{n}\right.$$

$$\left.\times \sum_{i=1}^n 1(s_{k-1} < x_i \leq s_k)\right\} dx_1\, dx_2 \cdots dx_n$$

$$= \lim_{m \to \infty} \frac{1}{\lambda^2} \sum_{1 \leq k < 2^m} s_k^{-2} \frac{\sigma^2}{n} \sum_{i=1}^n \int 1(s_{k-1} < x_i \leq s_k)\, dx_i$$

$$= \lim_{m \to \infty} \frac{1}{\lambda^2} \sigma^2 \sum_{1 \leq k < 2^m} s_k^{-2}(s_k - s_{k-1})$$

$$= \frac{\sigma^2}{\lambda^2} \int_\alpha^\beta \frac{1}{s^2}\, ds = (\alpha^{-1} - \beta^{-1})\frac{\sigma^2}{\lambda^2}. \qquad \square$$

**Acknowledgments.** The authors thank Professor Michael Kosorok for many discussions on the subject. They also thank the Editor Professor Bernard Silverman, the Associate Editor and two anonymous referees for useful comments and suggestions.

Y. Lan
100 Abbott Park Road
AP9A-1, R436
Abbott Park, Illinois, 60064
USA

M. Banerjee
G. Michailidis
Department of Statistics
University of Michigan
439 West Hall
1085 South University
Ann Arbor, Michigan 48109-1107
USA
E-mail: gmichail@umich.edu